\documentclass[11pt]{amsart}

\usepackage{pstricks}
\usepackage{amsmath}
\usepackage{amsfonts}
\usepackage{amssymb}
\usepackage{amsthm}
\usepackage{tikz}
\usepackage{tcolorbox}
\usepackage{pgf,tikz,pgfplots}
\usepackage{mathrsfs}
\usetikzlibrary{arrows}

\title{
Old and new results on the Furstenberg sets}

\author{Aihua Fan}
\address{
LAMFA, UMR 7352 CNRS, University of Picardie, 33 rue Saint Leu,80039 Amiens, France \& School of mathematics and Statistics, Central China Normal University,  430079 Wuhan, China }
\email{ai-hua.fan@u-picardie.fr}

\author{Herv\'e Queff\'elec}
\address{CNRS, Laboratoire Paul Painlev\'e, UMR 8524 \& Labex CEMPI (ANR-LABX-0007-01),
Universit\'e de Lille                          
Cit\'e Scientifique, B\^t. M2                                  
59655 Villeneuve d'Ascq Cedex,
FRANCE }
\email{herve.queffelec@univ-lille.fr} 

\author{Martine Queff\'elec}
\address{CNRS, Laboratoire Paul Painlev\'e, UMR 8524, \& Labex CEMPI (ANR-LABX-0007-01),
Universit\'e de Lille                          
Cit\'e Scientifique, B\^t. M2                                  
59655 Villeneuve d'Ascq Cedex,
FRANCE}
\email{martine.queffelec@univ-lille.fr}

\begin{document}

\begin{abstract}This paper is a complement to \cite{AHM}. It surveys the works on the Furstenberg set $S=\{2^{m}3^{n}: n\ge 0, m\ge 0\}$  and its random version $T$. 
We also present some new results. For example, it is proved that $T$  almost surely contains a subset of positive lower density which is $\frac{4}{3}$-Rider.  
It is also proved that a class of random sets of integers are Sidon sets when Bourgain's condition is not satisfied; this generalizes a result of Kahane-Katznelson.  
Some open questions about $S$ and $T$ are listed at the end of the paper.
\end{abstract}  


\maketitle
\newenvironment{demo}{\medbreak\begingroup\noindent{\bf Proof. }}{\hfill$\diamondsuit$\endgroup\goodbreak\medskip}

\newtheorem{Pb}{Problem}

\newtheorem{defi}{Definition}[section]
\newtheorem{theo}{Theorem}[section]
\newtheorem{prop}[theo]{Proposition}
\newtheorem{cor}[theo]{Corollary}
\newtheorem{lem}[theo]{Lemma}
\newtheorem{rem}{Remarque}
\newtheorem{ex}{Exemple}
\newtheorem{pb}{Problem}
\newcommand {\N}{{\mathbb{N}}} 
\newcommand {\Z}{{\mathbb{Z}}}
\newcommand {\E}{{\mathbb{E}}}
\newcommand {\R}{{\mathbb{R}}} 
\newcommand {\Q}{{\mathbb{Q}}} 
\newcommand {\C}{{\mathbb{C}}} 
\newcommand {\D}{{\cal D}} 
\newcommand {\T}{{\mathbb{T}}}
\newcommand {\F}{{\cal{F}}}  
\newcommand {\G}{{\cal{G}}}  
\renewcommand {\S}{{\cal{S}}}  
\newcommand {\A}{{\cal{A}}}  
\newcommand {\J}{{\cal{J}}}  
\newcommand {\B}{{\bf{BAD}}}  
\renewcommand {\l}{{\cal{L}}} 
\renewcommand {\r}{{\cal{R}}} 
\renewcommand {\b}{{\cal{B}}} 
\renewcommand {\k}{{\cal{K}}} 
\newcommand {\q}{{\cal{Q}}} 
\renewcommand {\P}{{\mathbb{P}}} 
\newcommand {\p}{{\cal{P}}} 
\newcommand {\noi}{\noindent}
\newcommand {\eps}{\varepsilon}

This survey, enriched by some new results,  is devoted to the innocent-looking Furstenberg set of integers
$$S=\{2^j\times 3^k :  j , k\in \N_0\}=:\{s_1<s_2<\cdots <s_n<\cdots\},$$
 and  its random analogue $T$,  to which we come a little later in this paper.
The distribution in $\mathbb{N}$ of the set $S$  had already been studied by Ramanujan and Hardy \cite{Ha} at the beginning of the last century. 
More generally, if $\{q_1,\dots, q_s\}$ is a finite set of mutually prime numbers, we denote by $S(q_1,\dots, q_s)$ the multiplicative semigroup generated by the $q_j's$; when $\{q_1,\dots, q_s\}$ are the $s$ first prime numbers, we get what arithmeticians call the $s$-friable numbers \cite{FT}. 
In the extreme case of $s=1$, the set reduces to a canonical Hadamard set which enjoys  many nice properties in ergodic theory and harmonic analysis (as is well-known). It is natural to explore which properties  remain true when we switch to 
$S(q_1, \dots, q_s)$ with $s\ge2$ and to other less lacunary sets.

We particularly focus on the semigroup $S(2,3)=:S$, the set of $2$-friable numbers, also called Furstenberg set because of its link with the famous metric conjecture formulated by  Furstenberg in 1960's  \cite{Fu1}, that we recall in Subsection  \ref{ss:furst}.  This conjecture and other  questions originate from the attempt to explicit the following feeling: ``expansions to bases 2 and 3 look very different!". Note that 2 and 3 can be replaced by multiplicatively independent bases $a$ and $b$ (i.e. $\log a/\log b\notin\Q$).

The aim of this study is to explore dynamical as well as harmonic analysis properties of $S$ and of its random analogue $T$.
This paper is mainly based on our previous one \cite{AHM}.  
Also, we take the opportunity to  elaborate on some improvements and some new considerations that we omitted or missed in \cite{AHM}. So, to some extent, this work may appear as a complement to \cite{AHM}.

We will make use of the classical notations: we put $\mathbb{T}=\R/\Z$ identified with $[0,1)$,
$\|x\|:=d(x,\Z)$, $x=\{x\}+[x]$; $e_n(x)=e(nx):=e^{2i\pi nx}$ for $n\in\Z$ and $x\in\T$.  We write $m$ for the Lebesgue measure on $\T$. $|E|$ will denote  the cardinality of a finite set $E$. 
By $u\lesssim v$ we mean $u\le Cv$ for some constant $C>0$.
We also adopt the Hardy notations 
$\asymp$ and $\sim$, Landau's $o$ and $O$ \cite[p.7]{HaWr} and the English convention $\N=\{1,2,\ldots, \}$. We denote $a\wedge b$  the gcd of the positive integers $a,b$.
Recall that a {\it Hadamard set} $E=(\lambda_n)\subset\N$ is a set that satisfies  $\lambda_{n+1}/\lambda_n\ge r$ for some $r>1$ and for all $n\ge 1$.

\section{Historical introduction and motivations}

The non-lacunary behaviour of the Furstenberg set $S$ seems to have been known for a long time. Let us first mention Hardy and Ramanujan's contribution,  followed by major related results.

\subsection{Ramanujan and others.}
 In his first  letter addressed to Hardy in 1913
(cf. \cite{Ha}, Chapter V), Ramanujan asserted without explanation that
 $$\displaystyle|S_N|:=|S\cap[1, N]|=\frac{\log 2N\times \log 3N}{2\log2\log3} +{1\over2}\cdot
 $$
 Hardy understood this formula as an approximation and  stated that there is no evidence to show how accurate Ramanujan supposed it to be.
 Chapter V of \cite{Ha} is devoted to this lattice-type problem.  
In two papers published in 1921 and 1922, Hardy and Littlewood considered this question and  proved the following estimation
\begin{equation}\label{eq:Hardy} 
|S_N|=\frac{\log 2N\times \log 3N}{2\log2\log3}+o(\log N/\log \log N).
\end{equation}
Ostrowski obtained quite similar results by different methods.
The estimation \eqref{eq:Hardy} implies that
$$|S_N|=  \frac{1}{2\log2\log3} \log^2 N + \frac{\log 6}{2\log2\log3} \log N + o(\log N),$$
By taking $N=s_n$, we can easily deduce (cf. \cite{AHM})
\begin{equation}\label{s_n}
 s_n \sim  \exp(C\sqrt{n})/\sqrt6, \quad {\rm with}\  C=\sqrt{2\log2\log3}.
 \end{equation} 
 So we see that $S$ has an intermediate sparseness with rate of increase
 both subexponential 
and superpolynomial. In \cite{AHM}, we have proved a more precise estimate of the remainder term in \eqref{eq:Hardy} by using diophantine approximation and discrepancy.
We will come back to this
in Section \ref{Sect:AA}.

\subsection{The number $\alpha=\log2/\log3$.} 
An attractive feature in this study is   the role played by the number $\alpha:=\frac{\log2}{\log3}=0.631\cdots$ and its diophantine properties.
This transcendental number is diophantine (non-Liouville),  which means that for some $C>0$ and $\rho\ge1$ 
$$\Vert q\alpha\Vert\geq C q^{-\rho} \hbox{\ for all }\ q\in\N.$$ 
Such a parameter $\rho$ is an irrationality exponent of $\alpha$, which is  thus called {\it $\rho$-diophantine.} Concerning a good value of $\rho$ for $\alpha$, G. Rhin proved that $\rho\leq 7.616...$ by using Pad\'e approximants \cite{Rh}; this estimate has been improved by Wang and Wu to
$\rho\le 4.11633052....$  \cite{WaWu}. 

Tijdeman \cite{Ti1,Ti2}, for more general $S=S(q_1,\ldots, q_s)$, established the following inequalities where $A,B$ are positive constants depending only on $q_1,\ldots, q_s$: 
$$s_n/(\log s_n)^A\lesssim s_{n+1}-s_n\lesssim s_n/(\log s_n)^B.$$ 
His proof used the Baker's results on linear forms of logarithms.

For the special case of $S=S(2,3)$, we have obtained   an improved form of these inequalities with {explicit} $A$ and $B$ in terms of any irrationality exponent 
$\rho$ of $\alpha$ \cite{AHM} (see also Theorem \ref{S_N}). 
The estimation on $|S_N|$ will be useful in considering the question whether $S$ is a Khintchin set, one of  the dynamical aspects of $S$ (see Section \ref{Sect:DA}).

\smallskip
\subsection{Two combinatorial properties of $S$.}\, \
  
\noindent  $\bullet$ G. Rauzy observed that the increasing sequence of powers of 2 or 3 is arranged according to a sturmian law. More precisely,
 there are one or two powers of $2$ between two consecutive powers of $3$, and the sequence $u$ with values in the alphabet $\{1,2\}$ obtained by coding one power of $2$ by $1$ and two powers of $2$ by $2$ is 
$$u:=1\ 2\ 1\ 2\ 1\ 2\ 2\ 1\ 2\ 1\ 2\ 2\ 1\ 2\cdots,$$
which is the characteristic sturmian sequence associated to $\alpha$ (see \cite{Ar} page 143, for equivalent definitions).
Indeed, $u_n=1$ means that there exist $k<j$ such that 
$3^k<2^j<3^{k+1}$, with $ j+k=n.$
Hence, $k\log3<j\log2<(k+1)\log3$, in other words
 $0<j\alpha-k<1$;  since $k=n-j$, we get $0<j(1+\alpha)-n<1$ then $0<j-n\beta<\beta$, where $\beta:=1/(1+\alpha)$ with continued fraction expansion  $[0;1,\alpha_1,\alpha_2,\cdots]$. Finally,
 $$u_n=1\Longleftrightarrow \{n\beta\}\in ]1-\beta,1),$$
which is one of the equivalent descriptions of the characteristic sturmian sequence with slope $\alpha\in(0,1)$.

\smallskip
\noindent $\bullet$  The continued fraction expansion
  of $\frac{\log 2}{\log 3}$ is related to $S$. The following observations and more on the convergents to $\alpha$ can be found in \cite{BDH}. Here are the first terms of the sequence $S$:
$$
 \ \boxed{2, 3, 2^2} \ {2\times3}\ \ \boxed{2^3,3^2}\
\ {2^2\times3}, \ {2^4}, \ {2\times3^2}, \ {2^3\times3}\,\ 
\boxed{3^3,  2^5}\ \ {2^2\times3^2},
$$
$$ {2^4\times3}, \ {2\times3^3}, \ {2^6}, \ {2^3\times3^2}, \ {3^4},
\  {2^5\times3}, \ {2^2\times3^3}, \ {2^7}, {2^4\times3^2}, \ {2\times3^4},
$$
$$ {2^6\times3}, \ {2^3\times3^3} \ 
\boxed{{3^5},2^8} \ {2^5\times3^2}, \ {2^2\times3^4}, \ {2^7\times3}, \ {2^4\times3^3}, \ {2^9}, \ {2\times3^5},$$
$$ {2^6\times3^2}, \ {2^3\times3^4}, \ {3^6},
 \ {2^8\times3}, \ {2^5\times3^3}, \
{2^2\times3^5 }, \ {2^{10}}, \ {2^7\times 3^2}, \  {2^4\times3^4},
$$
$$ 2^9\times 3, \ {2\times3^6}, \ {2^6\times3^3},  \ {2^3\times3^5} \ \boxed{2^{11}, 3^7}\ \dots\  
\boxed{2^{19}, 3^{12} }\ \dots\ \boxed{3^{19}, 2^{30} }\ \cdots
$$ 
We call a {\it pure pair}, any pair of consecutive terms $(3^p,2^q)$ or $(2^q,3^p)$ in the sequence $S$; these pairs lead exactly to the {\it Farey approximations} $(p/q)$ to $\alpha$ since 
$|q\alpha-p|<|q'\alpha-p'|$ for any $q'< q$, 
so that, if $p/q$ now is a {\it convergent} to $\alpha$, $(3^p,2^q)$ or $(2^q,3^p)$ according to the parity is a pure pair in $S$.
We easily identify the  first convergents to $\alpha$ as $1, 1/2, 2/3, 5/8,12/19,\dots$ whence 
$$\alpha=[0;1\ 1\ 1\ 2\ 2\ \cdots].$$
In view of Subsection 1.2,  the open question whether $\alpha$ could be a badly approximable number (i.e. with bounded partial quotients) becomes relevant since then $\rho=1$ is admissible.

\smallskip
\subsection{Furstenberg conjecture and dynamical properties of $S$.} \label{ss:furst}
 The study of $S$ is {motivated} by a famous, still open, conjecture of Furstenberg,
namely: 
{\it A continuous (i.e. atomless) probability Borel measure $\mu$ on $\mathbb{T}$, $\times 2\times 3$-invariant, must be equal to the Lebesgue measure $m$. } 
In terms of the Fourier transform, the conjecture of Furstenberg is stated as follows
\begin{equation*}(C)\quad \widehat\mu(2n)=\widehat\mu(3n)=\widehat\mu(n)\ \hbox{ for every}\ n\in\Z\Longrightarrow
     \mu=m.
\end{equation*}

This conjecture is related to the dynamics of the semi-group $S$, e.g. the distribution of the orbits $(s_nx)$ for various $x\in\T$.  Let $\Lambda=(\lambda_n)\subset\Z$; recall that the $\Lambda$-orbit of $x$ is the sequence $(\lambda_nx)$ of $\T$. 
\begin{defi} We say that $x$ is $\Lambda$-normal if
the sequence $(\lambda_nx)$ is uniformly distributed modulo 1.
\end{defi}
Thanks to the H. Weyl's theorem this means that 
\begin{equation}\label{HW} {1\over N}\sum_{n=1}^N e(h\lambda_n x)\to 0 \ \hbox{for every integer}\ h\not=0.
\end{equation}
It is known that for any increasing sequence $E=(\lambda_n)\subset \N$, (\ref{HW}) holds for almost all $x\in\T$ (\cite{KN}) hence the notation:
\begin{defi} We denote by $W(\Lambda)$ the negligible set of non-$\Lambda$-normal numbers.
\end{defi}

When $\Lambda=(q^n)$ for some $q\ge2$, the $\Lambda$-orbit of $x$ is just the orbit of $x$ under the action of the $q$-shift $\sigma_q:=x\mapsto qx$ on $\T$ and we speak of $q$-normal numbers.
In this case, any orbit $O_q(x):=(\sigma_q^nx)$ can be described by the $q$-adic expansion of $x$. It is thus easy to construct uncountable sets of $x$ with non-dense orbit and continuous $q$-invariant probability measures singular with respect to the Lebesgue measure like Bernoulli measures. The set $W(q)$ of non-$q$-normal numbers has a Hausdorff dimension 1 (which holds  more generally  for any lacunary (Hadamard) sequence in place of $(q^n)$, (cf. \cite{ET, Fa}).  The other way round, by the Birkhoff theorem, almost every orbit is $\mu$-distributed for any $q$-invariant ergodic measure $\mu$. 

Furstenberg proved in \cite{Fu1} that $(s_nx)$ is dense for every $x\notin\Q$, a first notable difference between $S$ and the lacunary sequences. At the end of his paper, he claimed that $W(S)$ is uncountable by suggesting an appropriate family of Liouville numbers inside. 
This is a first element regarding his conjecture.

Actually, if $W(S)$ were reduced to $\Q$, Furstenberg's conjecture {would hold.} Indeed, if it fails,
there exists a continuous probability measure $\mu$ on $\T$, $\times 2\times 3$-invariant, with $\widehat{\mu}(a)\neq 0$ for some $a\in  \mathbb{Z^\ast}$.
Then we get a contradiction:
$$0\not=|\widehat{\mu}(a)|=\limsup\Big|\frac{1}{N}\sum_{n=1}^N\widehat{\mu}(a s_n)\Big|\le\limsup\int_{\T}\Big|\frac{1}{N}\sum_{n=1}^N e(- a s_{n}x)\Big|d\mu(x)$$
$$\le\int_{\T}\limsup\Big|\frac{1}{N}\sum_{n=1}^N e(- a s_{n}x)\Big|d\mu(x)\le \mu(\mathbb{Q})=0.$$

Also in the 90' Bergelson asked whether $S$ could be a {\it recurrence set}, referring to the dynamical classification developed in \cite{Bo1, Ru}. We shall  prove in Section 4.2.3. that  this is not the case.

\smallskip
  Another dynamical classification is more adapted to the set $S$ and goes back to Hartman \cite{Har}. We briefly recall the definition that he introduced.
\begin{defi}  $E\subset\Z$ is a Ka-set if there exists a continuous measure $\mu$ on $\T$ such that
$\inf_{n\in E}|\widehat\mu(n)|>0$. (Here ``Ka" stands for R. Kaufman.)
\end{defi}
 Hartman observed that $W(E)$ is uncountable for any Ka-set $E:=(n_k)$. Indeed, if not,  we would have
$\lim_N{1\over N}\sum_{k\le N} e(n_kt)\to 0$ for every $t$ outside a countable set. 
As above, this leads to
\begin{equation}\label{bof}\lim_N{1\over N}\sum_{k\le N} \widehat\nu(n_k)\to 0
\end{equation} for every continuous measure $\nu$.
However, a continuous measure $\mu$ exists with $|\widehat\mu(n_k)|>\delta$ since $E$ is a Ka-set; the continuous measure $\nu=\mu\ast\tilde\mu$ ($\tilde\mu(A):=\overline{\mu(-A)}$) 
in turn satisfies $\widehat\nu(n_k)>\delta^2$ and (\ref{bof}), whence a contradiction.

Later R. Lyons \cite{Ly1}, investigating partial answers to the conjecture, asked whether $S$ could be a Ka-set (conjecturing the opposite), since this would be implied by the disproof of the conjecture. Badea and  Grivaux \cite {BG} gave a positive answer, recovering the uncountability of $W(S)$.

It becomes clear that the size or the shape of $W(S)$ are relevant in view to the conjecture. In \cite{AHM}, keeping in mind the Hadamard case, we prove that $W(S)$ is rather big since it has a Hausdorff dimension $>0.45$ (this has been recently improved to $\dim_{H}(W(S))=1$ by S.~Usuki \cite{Us}) and rather well distributed by constructing a {\it Rajchman measure} supported on $W(S)$.

\smallskip
\noi{\bf 1.5. Khintchin class.}  Recall that $(\lambda_n x)$ is almost surely uniformly distributed for any increasing sequence $E=(\lambda_n)\subset \N$. 
 As a consequence, for every {Riemann-integrable} function $f$ we have
\begin{equation} \label{Riemann0} \frac{1}{N}\sum_{n=1}^N f(\lambda_{n} x)\to \int_{\T} f dm \quad a.e.
\end{equation}
Khinchin conjectured that this holds for $L^\infty$-functions \cite{Kh}. 
Marstrand \cite{Ma}, refuting this conjecture, 
proved that 
  (\ref{Riemann0}) fails for $E=\N$, by  taking  $f={\bf 1}_O \in L^{\infty}(\T)$ with some well-chosen open set $O$. 
This led us  to coin the {{\it Khinchin class}} of  $E=(\lambda_n) \subset \N$ as 
\begin{equation}\label{khikhi} \mathcal{K}_E:=\Big\{f\in L^{1}(\T) :  \frac{1}{N}\sum_{n=1}^N f(\lambda_{n} x)\to \int_{\T} fdm \hbox{\quad a.e.} \Big\}
\end{equation}
 
 For the Furstenberg set $S$,  Marstrand  also proved in \cite{Ma}  that $\mathcal{K}_S$ contains $L^{\infty}(\mathbb{T})$ and, later, Nair \cite{Na} got 
$\mathcal{K}_S=L^{1}(\mathbb{T})$. 
 
 In \cite{AHM} we gave  an elementary proof of the inclusion $\mathcal{K}_S\supset L \log^{+}L$. The slight restriction on integrability is due to the following: we need some maximal function to be integrable, contrary to a fake result in Man\'e's book \cite{Mane} (Corollary 1.6, page 96), for which we will present a simple counterexample (Theorem \ref{jeje}). Actually,
using the classical Birkhoff {ergodic theorem,}  we proved a result which recovers Marstrand's result, nearly Nair's result and extends to $\times2\times3$-invariant probability measures 
(see Theorem \ref{thm:FQQ-K}). This highlights the fact that continuous and singular $\times2\times3$-invariant measures (counterexamples to the conjecture of Furstenberg), if any, must be carried by $W(S)$.
  
\smallskip
\noi{\bf 1.6. Lacunarity in  harmonic analysis.} Thin sets of integers play a fundamental role as \textit{spectrum} in harmonic analysis.
Here is a generic definition. If $X\subset  L^{1}(\T)$ is a Banach space and $E\subset \mathbb{Z}$ a subset of integers, we set 
$$X_E=\{f\in X : \widehat{f}(n)=0 \hbox{\ for}\ n\notin E\}.$$
It appears that functions in $X_E$, i.e. functions in $X$ with Fourier spectrum in a thin set $E$, enjoy better properties than a generic function in $X$ (we have  $X_E=\{0\}$ in  the extreme case $E=\emptyset$).  Here are some  special examples of thinness:  
\begin{itemize}
\item[(1)]  {\it Every $f\in L^{2}_{E}$ satisfies  $\Vert f\Vert_q \leq C_q \Vert f\Vert_2$ for some $2<q<\infty$.} 
We then say that $E$ is a $\Lambda(q)$-set 
(in short $L^{2}_{E}\subset L^q$).
\item[(2)]  {\it Every $f\in L^{\infty}_{E}$ satisfies $\sum_{n}|\widehat{f}(n)|<\infty$.}
We then say that $E$ is a Sidon set (in short $L^{\infty}_{E}\subset A(\T)$, the Wiener algebra).
 \item[(3)]
 {\it Every $f\in L^{\infty}_{E}$ satisfies $\sum_{n}|\widehat{f}(n)|^p<\infty$ for  some fixed $1\leq p<2$}. We then say that $E$ is a $p$-Sidon set ($L^{\infty}_{E}\subset A_p(\T)$). 
\end{itemize}

It is well known that Hadamard lacunary sets enjoy all these properties (for $p=1$ and for all $q<\infty$). The Furstenberg set $E=S$, which is less lacunary, is an interesting candidate to be considered  from this harmonic analysis viewpoint. 
Such properties for $S$ seem {poorly known}, contrary to the well-understood sumset $\{2^j+3^k\}$ \cite{ LoR} (we could say that the characters on $\T$ digest better the sum of integers than the product!).  
However, in 1970's,  Gundy and Varopoulos \cite{GV} proved  that 
{\it the set $S$ is  $\Lambda(q)$ for all $q>2$}
(and even more). 
We have produced a proof in \cite{AHM} when $S=S(q_1,\dots,q_s)$ with a careful analysis of the constants and their dependence on the parameter $s$. The proof uses the {square function} of a martingale and Burkholder's inequalities. 

In a fascinating way, combinatorics reappears in this study of thin sets $E$; in particular, the rate of growth of the cardinality $|E_N|:=|E\cap[1,N]|$ reveals a main tool \cite{AHM}, which 
immediately implies that {\it $S$ is not a Sidon set, even not $p$-Sidon for $p<4/3$}; we will come back to this in the following. Another measuring tool, less easy to implement, consists in enumerating the relationships in $E$: $\sum\epsilon_jn_j=0$ with $\epsilon_j\in\{0\pm1\}$ if $E:=(n_j)$. This is due to the use of Riesz products which we cannot do without. When $E=S$, once again, we are faced with confrontational additive and multiplicative properties.

\smallskip
\noi{\bf 1.7. Random version of Furstenberg set.}
 In order to understand which properties of $S$ depend on its arithmetic structure and which on its sparseness, we parallelly considered  in \cite{AHM} a random version $T$ of $S$: 
 let $(\delta_k)$ be a sequence of numbers with $0\le \delta_k<1$ and $\sum \delta_k =\infty$, and let $(\xi_k)_{k\geq 1}$ be a sequence of $(0$-$1)$-valued independent random variables with {\it expectation}
 $\E(\xi_k)=:\delta_k$; put
 $$R=R(\omega)=\{ k: \xi_{k}(\omega)=1\}.$$
 This is a random set of integers.
 The random version $T$ of the Furstenberg set corresponds to the special choice 
  $\delta_k={\log k}/{k}$.
  Indeed, almost surely $T$ has the same growth as  $S$ : with $T_N:=T\cap [1,N]$,
 $$|T_N|\sim\E(|T_N|)=\sum_{k=1}^N \frac{\log k}{k}\sim \int_{1}^N \frac{\log t}{t}dt=\frac{\log^2 N}{2},$$
 where the first ``$\sim$'' is an unconventional   law of large numbers (the normalizer is  $\E(|T_N|$ which is not proportional to $n$). See \cite{BF2005} for a proof of this law of large numbers. 
 
 Now we have lost the arithmetical property of the initial set $S$, whence the question of the remaining, or additional, properties of $T$.  
 Here are some results.
 
1.  $T$ is almost surely Hartman-distributed ($W(T)=\Q$), see Bourgain \cite{Bo2}. More generally, this holds for the random set $R$ as soon as $\delta_k \downarrow $ and  $m_N/  
 \log N\to \infty$ where $m_N=\sum_{k=1}^N \delta_k$. It is pointed out in \cite{FS} that the monotonicity of $\delta_k$ can be weakened to 
 $\sum_{n=1}^N |\delta_n-\delta_{n+1}|=o(m_N)$. In  this work, we  prove a converse of Bourgain's result (extending a result of Kahane and Katznelson):  if $\delta_k \downarrow $ and $m_N=O(\log N)$, then $R$ is almost surely Sidon, and hence not Hartman-distributed (see Theorem \ref{mieux que kk}).
 
2.
 $T$ is almost surely $p$-Rider for each $p>4/3$ \cite{AHM}. We were not able to prove the same result for $S$, but will prove here that $S$ is $3/2$-Rider
 (see Theorem \ref{gene} and Theorem \ref{decept}).
 
3.
$T$ contains almost surely a subset $T'$ of positive lower density which is $\Lambda(p)$ for each $2< p<\infty$ \cite{LQR2}. We conjecture that $T$ itself is $\Lambda(p)$ for each $p<\infty$. 

\bigskip 
 The following comparative table is very meaningful.
 
 \bigskip\noi
 {\begin{tabular}{|c|c|}\hline

{DETERMINISTIC SET $S$} \ \ \  &  { RANDOM SET $T$} \\
 \ \ \  &\\
 \hline

{ Weak lacunarity} of $S$\ \ \  &  {Similar lacunarity} for $T$\\
\ \ \  &\\
 \hline
   
{Hausdorff} dimension of $W(S)=1$ \ \ \   & $W(T)=\mathbb{Q}$  \\  
  \ \ \  &\\
\hline

$W(S)$ {supports} a Rajchman measure\ \ \  & $W(T)=\mathbb{Q}$  \\ 
 \ \ \  &\\
\hline

 $S$ {is a Ka-set} \ \ \   & $T$ is not a Ka-set\\ 
  \ \ \  &\\
  
  \hline

 $S$ {is a Bohr-closed} \ \ \   & $T$ is {Bohr-dense}\\ 
  \ \ \  &\\

\hline

$\mathcal{K}_S=L^{1}$ \ \ \  &  ``every'' $f\in\mathcal{K}_\N$ is a.s. contained in $\mathcal{K}_T$\\ 
   \ \ \  &\\
 \hline

 $S$ {is $\Lambda(p)$} \ \ \  &  {A big subset  of $T$ is $\Lambda(p)$.}\\ 
  \ \ \  &\\
 \hline 
 
$S$ is  p-Rider for  $p=3/2$, \ \ \   &  {$T$ is} $p$-Rider for $p>4/3$, \\  
  \ \ \  &\\ 
  how about $4/3<p<3/2$ ? \ \ \   &  but not $p$-Rider for $p<4/3$. \\  
  \ \ \  &\\ 
\hline
\end{tabular} }

\section{Arithmetical aspects}\label{Sect:AA}
It is easy to prove that $|S_N|=:|S\cap [1,N]|\asymp (\log N)^2$ whence the rate of growth
 $s_n\asymp \exp(\sqrt{n})$. 
Going further needs sharp estimates and more work. More generally, 
for $S(q_1,\dots, q_s)$, the multiplicative semi-group generated by coprime integers $q_1,\ldots, q_s$,  Marstrand \cite[p.545]{Ma}  gave the following  asymptotic estimate 
  \begin{equation}\label{eq:Marstrand-Est}
 |S(q_1,\dots, q_s)\cap[1,n]| = K_s (\log n)^s  +O((\log n)^{s-1}) \quad
 {\rm as } \ \  n\to\infty,
   \end{equation}
 where 
 $$K_s = \frac{1}{s!} \prod_{j=1}^s\frac{1}{\log q_j}
$$
which  is the volume of the simplex  $\{x=(x_j)\in\R^s :\  \sum_{1}^s\lambda_j x_j\le1,\  x_j\ge0\}$, with here $\lambda_j=\log q_j$.

  Our first task in \cite{AHM}, in the case of $S$, consisted in an improved remainder term with respect to Hardy's one, namely $o(\log N/\log \log N$), by using diophantine approximation properties 
of the  number $\alpha=\frac{\log 2}{\log 3}$ and the discrepancy of the sequence $(n\alpha)$.
  
  According to Gelfond  \cite{Ge}, 
$\frac{\log a}{\log b}$ is either rational ($a=N^p, b=N^q$) or transcendental, so that $\alpha$ is indeed transcendental. But it is {not Liouville} meaning that 
there exists $\delta>0$ and  $\rho<\infty$ such that
$$\forall q\in \mathbb{N}, \qquad  \Vert q\alpha\Vert\geq \delta q^{-\rho} $$
Such a parameter $\rho\geq 1$ is called {\it an irrationality exponent} of $\alpha$, which is  then called  {$\rho$-diophantine.} Concerning a good $\rho$ for $\alpha$, as said in the introduction, G. Rhin \cite{Rh} obtained
$\rho\leq 7.616...$ and many years later Wang and Wu \cite{WaWu} got the improved estimate $\rho\leq 4.11633052....$.
Based on this estimate, the following result is proved.

 \begin{theo}[\cite{AHM}]\label{S_N} Assume that $\alpha$ is $\rho$-diophantine; let $\delta:=\frac{\rho}{\rho +1}$ . Then 
\begin{equation}\label{eq:FQQ}
|S_N|= {1\over 2\log2\log3}\log^2N + 
 {\log6\over 2\log2\log3}\log N +O((\log N)^\delta).  
 \end{equation}
 Or else
 $$|S_N|=\frac{\log 2N\times \log 3N}{2\log2\log3}+O((\log N)^\delta).$$
 We can take $\rho= 4.11633052....  $ and $\delta:= 0.80454645\cdots$
\end{theo}


\bigskip
In our case $S=S(2,3)$, we could improve the inequalities of Tijdeman for $S(q_1,\ldots, q_s)$,  with here {explicit}  constants depending once more on an irrationality exponent 
$\rho$ of $\alpha$.

\begin{theo}[\cite{AHM}] \label{ss} We have that \textnormal{(}with $2r=1/(\rho+1)$\textnormal{)}
\begin{equation}\label{one}{1\over(\log s_n)^{\rho}}\lesssim{s_{n+1}-s_n\over s_n}\lesssim{1\over n^{r}}\lesssim{1\over(\log s_n)^{2r}}.
\end{equation}
 We can take  $\rho\sim4.116$ and  $2r\sim0.1954$.
 In particular, $s_{n+1}-s_n \to \infty$ and  $s_{n+1}/s_n \to 1$.
\end{theo}

\noi{\bf Remarks.} 1. To what level can the right inequality in \eqref{one} be reinforced was a question raised by Tijdeman who observed that $\log s_n$ is best possible. We can add a small precision to (\ref{one}): for infinitely many pairs $(s_n,s_{n+1})$, related to the approximations to $\alpha$, we have
$${s_{n+1}-s_n\over s_n}\le {2\log2\log3\over\log s_n}\cdot$$
This is due to the best approximation property of the convergents \cite{AHM}. 
\\
2. It then remains to locate those pairs. If $\alpha$ is a badly approximable number, thoses pairs appear with bounded gaps and the inequalities can be refined.  If not, the authors in \cite{BDH} observe that $S$ contains arbitrarily long intervals in a {\it geometric} progression: {\it For infinitely many $n$ there exist $a_n:=2^{-q_n}3^{p_n}\in\Q$ and $j:=j_n\ge1$ such that $s_{k+1}/s_k=a_n$ for $j\le k\le j+n-1$}.
 
\section{Dynamical aspects}\label{Sect:DA}

We now turn to the dynamics of the sequence $S$ and the description of the $S$-orbits, i.e. of $\{s_nx\}$ for $x\in\T$. Finite orbits do exist and can be described.
A topological study has been achieved by Furstenberg who proved the following rigidity result.
\begin{theo} [{\cite{Fu1}}]
An infinite closed subset of $\T$ which is $S$-invariant must be $\T$ itself.
\end{theo}


In particular, every infinite $S$-orbit is dense in $\T$. This very interesting result reveals a deep difference between the $S$-action and the rich dynamics of $q$-shifts ($q\ge2$) where infinite non-dense orbits exist.
To go further, {we need a more elaborate  analysis of the distribution of the points $s_nx$.}  
We say that the sequence of numbers $(u_k)\subset\T$ is uniformly distributed if, for any $0\le a<b<1$,
$${1\over N}|\{1\le k\le N,\ u_k\in[a,b)\}| \to b-a;$$
a very useful criterion due to H. Weyl asserts that $(u_k)$ is uniformly distributed if and only if, for every $h\not=0$,
$${1\over N}\sum_{1\le k\le N}e(hu_k)\to 0.$$
We shall also make use of the following definition.
\begin{defi} A sequence $E:=(n_k)$ of integers is Hartman-distributed if  
$${1\over N}\sum_{1\le k\le N}e(n_kx)\to 0\ \ \forall x\in\T,\ x\not=0.$$
\end{defi}
The negligible set $W(E)$ of $x\in \T$ such that $(n_kx)$ is not uniformly distributed comes into play. Observe that
``$E$ is Hartman-distributed" means ``$W(E)=\Q$".

\subsection{Khintchin class}
Consider an arbitrary increasing sequence of integers $E=\{\lambda_n\}_{n\geq 1} \subset \mathbb{N}$. It is well known that  $(\lambda_n x)$ is a.e. uniformly distributed. 
 Consequently, for every {Riemann-integrable} function $f$ we have
\begin{equation} \label{Riemann} \frac{1}{N}\sum_{n=1}^N f(\lambda_{n} x)\to \int_{\mathbb{T}} fdm \quad a.e. 
\end{equation}
Marstrand \cite{Ma} 
proved that when $E=\{\lambda_n\}=\mathbb{N}$, there are bounded functions $f\in L^{\infty}(\mathbb{T})$
 such that (\ref{Riemann}) fails. We define 
  the {{\it Khinchin class}} of  a set $E=\{\lambda_n\} \subset \mathbb{N}$ as that of those Lebesgue integrable functions such that (\ref{Riemann}) holds.

 Marstrand \cite{Ma}  proved that 
 \begin{equation}\label{mama}\mathcal{K}_S \supset L^{\infty}(\mathbb{T})\end{equation} 
 and later Nair \cite{Na}  proved that 
$\mathcal{K}_S=L^{1}(\mathbb{T})$.  Thus the Khintchin class $\mathcal{K}_S$ is completely determined,
but the determination of $\mathcal{K}_\mathbb{N}$ is not complete. 
 {Koksma} \cite{Ko} proved the following criterion for   $E=\mathbb{N}$:  
 $$\sum_{|k|\geq 3} |\widehat{f}(k)|^2 \big(\log\log |k|)^{3}<\infty\Longrightarrow f\in \mathcal{K}_{\mathbb{N}}.$$  
 
 For $E=\{q^n\}$ with an integer $q\ge 2$, we have $\mathcal{K}_{E}=L^1(\mathbb{T})$ by the Birkhoff ergodic theorem. 
 But for a general lacunary set $E$, the determination of $\mathcal{K}_E$ is unknown.
 In many cases, lacunary sequences  share nice properties. 
 The sequence $\{\lambda_n\}=\{2^{2^n}\}$ is very lacunary, but J. Rosenblatt \cite{Ros} proved that there exists $f\in L^\infty$ such that (\ref{Riemann}) fails, in other words
 $\mathcal{K}_{\{2^{2^n}\}}\not\supset L^\infty$. 
 For any Hadamard lacunary sequence $\Lambda=\{\lambda_n\}$, a result of Cuny and Fan \cite[Theorem C, p.~2728]{CF}  implies  that 
 $$
      \sum_{n>N}|\widehat{f}(n)|^2 = O\Big(\frac{1}{\log^{1+\epsilon} N}\Big)
       \Longleftrightarrow \omega_{f,2}(\delta) = O\Big(\frac{1}{|\log \delta|^{1/2 +\epsilon/2}}\Big) \Longrightarrow f \in \mathcal{K}_\Lambda
 $$
 where $ \omega_{f,2}(t) $ denotes the modulus of  $L^2$-continuity of $f$ defined by 
 $$
            \omega_{f,2}(t) =\sup_{|t|\le \delta} \|f(\cdot) - f(\cdot -t)\|_2.
 $$
 Actually, the above condition on $\omega_{f, 2}(\delta)$ (with assumption $\int f =0$) implies that $\{f(\lambda_n x)\}$ is a convergence system, meaning that
 $$
        \sum_{n=1}^\infty |a_n|^2 <\infty \Longrightarrow \sum_{n=1}^\infty a_n f(\lambda_n x) \ \ {\rm converges \ a.e.}
 $$
 and the exponent $1/2$ is the best possible \cite[Proposition 5.2]{CF}.
 
Nair's proof of $\mathcal{K}_S=L^{1}(\mathbb{T})$ is based on an ergodic theorem for amenable group actions due to Bewley \cite{Bewley1971}.
Using the classical ergodic theorem,  it is possible to give  a {simple proof} of Marstrand's result $\mathcal{K}_S \supset L^{\infty}(\mathbb{T})$, and even
 nearly of Nair's result, with an extension to a class of $\times2\times3$-invariant probability measures. 
 The result is stated as follows.

\begin{theo}[\cite{AHM}] \label{thm:FQQ-K} Let $\mu$ be a $\times2\times3$ invariant probability measure, {ergodic} for one of both shifts. If $f\in L\log^+L$, then, 
$$A_Nf(x):={1\over N}\sum_{n\le N}f(s_nx)\to \int fd\mu\quad   \mu{\rm-}a.e.$$
\end{theo}
As a consequence, if a {``Furstenberg-exotic''} measure $\mu$ exists, namely continuous $S$-invariant and $\mu\neq m$, then it must be supported by $W(S)$ (actually $\mu$ is of zero dimension according to a result of Rudolph \cite{Rud}):  indeed, we can suppose that in addition $\mu$ is $S$-ergodic by invoking the $S$-ergodic decomposition; now choose $a\not=0$ such that $\widehat\mu(a)\not=0$ and apply theorem \ref{thm:FQQ-K} with $f=e_a$:
${1\over N}\sum_{n\le N}e(as_nx)\to\widehat\mu(a)$ $\mu$-a.e. and such $x's$ belong to $W(S)$.

\smallskip
A key step in the proof of Theorem \ref{thm:FQQ-K} is the following generalization of Birkhoff's ergodic theorem.
 
 \begin{theo}[\cite{AHM}]\label{thm:GBirkhoff}
 Let $(X, \mathcal{B}, \mu, T)$ be a measure-preserving dynamical system.  Let $(f_n)$  be a sequence of integrable functions. Suppose that 
 \\
 \indent  {\rm (1)} \  $f_n\to 0$ a.e. and  $\| f_n\|_1\to 0$;
 {\rm (2)}\  $\sup_{n} |f_n|$ is integrable.\\
Then almost everywhere we have  
    $$ 
    \frac{1}{N} \sum_{n=0}^{N-1} f_{n}(T^{N-n}x) \to 0.$$ 
    \end{theo}
    
 Ma\~{n}\'e claimed the same conclusion in \cite[pages 96--97]{Mane}  without assuming that $\sup_{n} |f_n| $ is integrable.  His proof   presents a gap. We show here that this extra assumption is mandatory in the general case and cannot be dropped.
We  detail for that a counterexample due to F.~Rodriguez-Hertz (kindly indicated to us by  L.~Flaminio \cite{Fl}).  The idea lurking behind is that of ``shrinking target". We will give a self-contained and elementary proof.
\begin{theo}\label{jeje} Let $(X, \mathcal{B}(X), \mu, T)$ be the Bernoulli system where $X=\{0,1\}^{\N}$, $\mu$ is the  symmetric Bernoulli measure  and $T:X\to X$  the left shift defined by  $T((x_i))=(x_{i+1})$. 
     Let $(\lambda_n)$ be a non-decreasing sequence of positive integers and 
     let {$E_n=\{x: x_1=\cdots=x_{\lambda_n}=0\}$,} a cylinder of length $\lambda_n$. Suppose that 
   \begin{displaymath} 
   \sum \frac{1}{2^{\lambda_n}}=\infty, \qquad \frac{n}{2^{\lambda_n}}\to 0.
   \end{displaymath} 
   (one choice is $\lambda_n=[\log(n\log n)/\log2]$).
  The sequence of functions $(f_n)$  defined by $f_n=n1_{E_n}$ satisfies the condition (1) in Theorem \ref{thm:GBirkhoff}, but 
   $$a.s \quad \varlimsup_{N\to \infty} \frac{1}{2N} \sum_{n=0}^{2N-1} f_{n}(T^{2N-n}x)\geq 1/2. $$
      \end{theo}
      \begin{proof} We first observe that $E_n \downarrow \{0^\infty\}$. It follows that $f_{n}(x)\to 0$ for all $x\not=0^\infty$. On the other hand, 
       $\Vert f_n\Vert_1=n\mu(E_n)= {n}/{2^{\lambda_n}} \to 0.$
      Thus the condition (1) in Theorem \ref{thm:GBirkhoff} is met.
      
       In order to prove the announced inequality, we shall use the following form of the Borel-Cantelli lemma \cite[p.368]{Re}.
      Let $(A_n)$ be a sequence of events. We have  $\mu(\varlimsup A_n)=1$ if the following conditions are satisfied:
     \begin{equation}\label{eq:BC} \sum_{k=1}^\infty \mu(A_k)=\infty, \qquad 
     \varliminf_{n\to \infty} \frac{\sum_{1\leq k, l\leq n} \mu(A_k \cap A_l)}{(\sum_{1\leq k\leq n} \mu(A_k) \big)^2}=1.
     \end{equation}
       Observe that  
       $$\frac{1}{2N} \sum_{n=0}^{2N-1} f_{n}(T^{2N-n}x)\geq \frac{1}{2N} f_{N}(T^{N}x)= \frac{1}{2N} N 1_{E_N} (T^{N} x)= \frac{1}{2}  1_{A_N} (x) $$ where
  { $$A_N=\{x: x_{N+1}=\cdots= x_{N+\lambda_N}=0\} = T^{-N}(E_N).$$}
      We only need to check that $A_N$'s satisfy the condition (\ref{eq:BC}) to validate the Borel-Cantelli lemma. First, $\sum \mu(A_n)=\sum \mu(E_n)=\sum {2^{-\lambda_n}}=\infty$. Second,  set 
      $$J_n=\{ (k,l): 1\le k, l\le n ,  \   \mu(A_k\cap A_l)\not=\mu(A_k) \mu(A_l)\}.$$
Since  $ \mu(A_k\cap A_l)=\mu(A_k) \mu(A_l)$ for $(k,l)\notin J_n$, we get
\begin{equation}\label{eq:BC1}
   \sum_{1\le k,l\le n}  \mu(A_k\cap A_l) 
   =  \sum_{(k, l)\in J_n} [\mu(A_k\cap A_l)-\mu(A_k)\mu(A_l)] + \sum_{1\le k, l\le n}\mu(A_k)\mu(A_l).
\end{equation}
       Now assume $1\le k,l\le n$ and $l\ge k$ (the case of $k\ge l$ can be similarly dealt with). We  distinguish 
       two cases:\\
        \indent (i) \ if $l>k+\lambda_k$, then $(k,l)\in J_n$ because $A_k$ and $A_l$ are independent and hence
        $\mu(A_k\cap A_l)= \mu(A_k)\mu(A_l);$\\
       \indent (ii)   if $l\leq k+\lambda_k$, we have 
    {$A_k\cap A_l=\{ x_{k+1}=\cdots = x_{l+\lambda_l}=0\}$} so that 
{ \begin{displaymath} \mu(A_k\cap A_l) =\frac{1}{2^{l-k}} \cdot \frac{1}{2^{\lambda_l}}\leq \frac{1}{2^{l-k}} \mu(A_k).\end{displaymath} }
      It follows from (i) and (ii) that 
      \begin{equation}\label{eq:BC2}
      \sum_{(k,l)\in J_n} \mu(A_k\cap A_l)\leq 2 \sum_{1\leq k \leq n,\atop k\leq l\leq k+\lambda_k}  \frac{1}{2^{l-k}}\mu(A_k)\leq 4 \sum_{k=1}^n \mu(A_k)= o\big(\sum_{k=1}^n \mu(A_k)\big)^2.
      \end{equation}
      Similarly, we have 
      $$\sum_{(k,l)\in J_n} \mu(A_k)\mu(A_l) = 2 \sum_{1\leq k \leq n,\atop k\leq l\leq k+\lambda_k}  \mu(A_k) \mu(A_l)
      \le 2\sum_{1\le k\le n} \mu(A_k) \frac{\lambda_k}{2^{\lambda_k}}$$
      \begin{equation}\label{eq:BC3}
      =o\big(\sum_{1\le k\le n} \mu(A_k)\big). 
      \end{equation}
     The second condition in \eqref{eq:BC} is thus implied by \eqref{eq:BC1}, \eqref{eq:BC2} and \eqref{eq:BC3} and the fact {that}
     $\sum \mu(A_k)=\infty$. 
       \end{proof}

\subsection{Dimension of $W(S)$ and Rajchman measure on $W(S)$}
The previous remarks motivate a closer examination of $W(S)$.
 We can prove: 
 \begin{theo}[\cite{AHM},\cite{Us}]\label{baddis} The set $W(S)$ satisfies:
 \begin{enumerate}
 \item $\dim_{H} (W(S))   =1$.
 \item $W(S)$  supports a \textnormal{Rajchman} {probability} measure $\mu$, more explicitely,\\
 $ \widehat{\mu}(h)=O(1/\log\log |h|)$ as $h\to\infty$.
 \end{enumerate}
\end{theo}

In \cite{AHM}, we have only proved $\dim_{H} (W(S))   \ge 0.451621$, but this estimate has recently been improved to $\dim_{H} W(S)=1$ by Usuki \cite{Us}, using a method similar to ours.
The second assertion of Theorem \ref{baddis} indicates that $W(S)$ is {not that porous.} So is the (uncountable) set $\mathcal{L}$ of Liouville numbers, though $\dim_{H} (\mathcal{L})=0$ (\cite{Bl}). On the opposite, the Cantor ternary set $K$  (satisfying  $\dim_{H} (K)=\frac{\log 2}{\log 3}$) is {porous:} it supports NO Rajchman measure (Kahane-Salem \cite{KS}). 

 Now we are able to precise the $M_0$-property of the set $W(S)$ by showing that we can hardly do better than the decay stated in Theorem \ref{baddis} .
 
 \begin{theo}\label{baddis2} 
 For any  probability measure $\mu$  supported on $W(S)$,  $\widehat{\mu}(h)$ cannot decay as $1/(\log\log h)^\alpha$ as $h\to \infty$ as soon as $\alpha>1$.
 \end{theo}
 
    We are going to prove  Theorem \ref{baddis2} basing ourselves  on the classical criterion of Davenport-Erd\" os-LeVeque below. 
  \begin{prop} [cf \cite{Bug 2}, Lemma 1.8.] \label{del} Let $(\Omega, \mathcal{A}, \mu)$ be a probability space and $(X_n)$ a sequence of complex valued random variables with $|X_n|\leq 1$.  Then 
  the averages $Z_n:=\frac{1}{n}\sum_{j=1}^n X_j$ converge $\mu$-a.e. to $0$ under the assumption
 $$\sum_{n=1}^\infty \frac{1}{n} \int_{\Omega} |Z_n|^2 d\mu<\infty. 
 $$
   \end{prop} 
   
 \noindent   {\it Proof of Theorem \ref{baddis2}.}
   We will show that if $\widehat{\mu}(h)=O\big(1/(\log\log h)^\alpha\big)$ for some $\alpha>1$, then $\mu(W(S))=0$,  so that $\mu$ cannot be supported on $W(S)$.
   To this end, we  only need to show that  if $\ell$ is a non-zero integer, the averages $A_{n}^{(\ell)}(x):=\frac{1}{n}\sum_{j=1}^n e(\ell s_j x)$ 
   tend almost everywhere to $0$, which implies that $(s_nx)$ is uniformly distributed $\mu$-a.e.  
   
   We can assume that $\ell=1$ (the proof is the same for all $\ell$'s). Denote simply $A_n^{(\ell)}$ by $A_n$. We are going to apply 
   Proposition \ref{del} by just checking  its assumption for $X_{n}(x)=e(s_n x)$. Clearly, by expanding $|A_n|^2$ we get
   $$ \int_{\T} |A_n|^2 d\mu=\frac{1}{n^2} \Big(n+\sum_{ j\neq k, 1\leq j,k\leq n} \widehat{\mu}(s_k-s_j)\Big)
   \le\frac{1}{n^2} \Big(n+2\sum_{ 1\leq j<k\leq n} |\widehat{\mu}(s_k-s_j)|\Big).$$
  But we know that, for $j<k$ and some numerical constant $A$, it holds
  $s_k-s_j\geq s_k -s_{k-1} \gtrsim s_k/(\log s_k)^{A}$ (cf. Theorem \ref{ss}),
  so that 
 $$ \log\log (s_k-s_j)\gtrsim \log\log s_k\gtrsim \log (\sqrt {k})\gtrsim \log k$$
 and hence
  $ |\widehat{\mu}(s_k-s_j)|\lesssim 1/ (\log k)^\alpha.$
  This gives us 
  $$ \int_{\T} |A_n|^2 d\mu\lesssim \frac{1}{n^2} \big(n+\sum_{ 1\leq j<k\leq n}1/ (\log k)^\alpha\big)\lesssim \frac{1}{n}+ \frac{1}{n} \sum_{k=2}^n 1/ (\log k)^\alpha\lesssim \frac{1}{(\log n)^\alpha}.$$
 So, the assumption of  Proposition \ref{del} is met.
 \hfill$\Box$

 Theorem \ref{baddis2} is a manifestation of an uncertainty principle for measures on $\T$: the support and the spectrum of $\mu$ cannot be too small at the same time.

\section{Harmonic analysis of thin sets}\label{Sect:HA}


\subsection{Combinatorics  in harmonic analysis}
Combinatorics plays a fundamental role in harmonic analysis of thin sets, under various aspects such as independence, sparseness or  arithmetic relations. The cumulative function $N\to |E_N|:=|E\cap[1,N]|$ of a subset $E\subset\Z$ accounts for the sparseness of the set and provides necessary conditions for $E$ to enjoy one of the harmonic properties we are interested in, as recalled below.

\begin{defi} Let $E$ be a set of positive integers, and $X$ a Banach space of integrable functions on the circle (e.g. $X=L^p, 1< p<2
$). We say that $E$ is $X$-Paley if the Fourier transform of any  function $f\in X$, once restricted to $E$, is square-summable, i.e.
$$ \widehat{f}_{|E} \in \ell^2.$$
Then, there is a smallest constant $P(X,E)>0$, the Paley constant of the pair $(X,E)$, such that  the following Paley-type inequality holds:
$$\Vert \widehat{f}_{|E} \Vert_2\leq P(X,E) \Vert f\Vert_X,\quad  \forall f\in X.$$ 

\end{defi}

\begin{prop}[\cite{Ru2}] \label{maille}
The following holds.

\indent {\rm 1)} $E$ is $H^1$-Paley if and only if  $|E\cap[N,2N]|=O(1)$. 

\indent {\rm 2)} If $E$ is $p$-Sidon, then $|E_N|\le c(\log N)^{p\over 2-p}$.

\indent {\rm 3)}  If $E$ is $\Lambda(p)$ for $p>2$, then $|E_N|\le 4 \lambda_{p}(E)^2 N^{2/p}$.
\end{prop}

As an easy corollary of 1) (this indeed motivated the definition),  the set $\{2^n\}$, more generally a Hadamard set,  is $H^1$-Paley. The known estimate of $|S_N|$, together with 2) and 3),  leads to the first negative results on the $S$.
\begin{prop} \, \  \\
\indent {\rm (i)} \  $S$ is not Sidon and even more, $S$ is not  $p$-Sidon for $p<4/3$.\\
\indent {\rm (ii)}\ $S$ is not $H^1$-Paley.
\end{prop}
\subsubsection{Quasi-independent and Pisier criterion for Sidon sets}
\begin{defi} A set $E\subset\Z\backslash\{0\}$ is said to be quasi-independent if, for all \textnormal{distinct} elements $x_1,\ldots,x_n\in E$ and for all $\eps_1,\ldots,\eps_n\in\{-1,0,1\}$, the \textnormal{relation} $\sum\eps_kx_k=0$ implies that all $\eps_k=0$. 
\end{defi}

The quantity $\sum_{k}|\varepsilon_k|$ is called the \textnormal{length} of the relation. A quasi-indepen\-dent set is hence a set which contains no relation of positive length. 



Quasi-independent sets are prototypes of Sidon sets.  Since the sidonicity is preserved by finite union (Drury's theorem \cite{Dr}), it is conjectured that a Sidon set could be a finite union of quasi-independent sets. In this direction, a breakthrough has been made by Pisier \cite{Pi} with the following characterization. 
\begin{theo}[\cite{Pi}]
{\it $E\subset\Z\backslash\{0\}$ is a Sidon set if and only if there exists $\delta>0$ such that, from every finite subset $A\subset E$,  a quasi-independent set $B$ can be extracted from $A$ with $|B|\ge\delta|A|$.} 
\end{theo}
An analogue for $p$-Sidon and $p$-Rider sets will appear and be used later.

\subsubsection{Spectrum of continuous measures}
The conjecture of Furstenberg, actually, raises questions on continuous measures on $\T$ and their spectrum. Continuous measures can be described just as well in terms of support (actually of annihilating set) as in terms of Fourier coefficients thanks to Wiener criterion (a scent of the uncertainty principle).
Russel Lyons (\cite{Ly0}) obtained the following new characterization:\\ 
{\it
A measure $\mu\in M(\T)$ is continuous if and only if there exists $(n_k)\subset\Z$, $|n_k|\to\infty$ such that
$\widehat\mu(n_k\ell)\to 0$ for every $\ell\in\Z^*$.
}

But nothing can be said in general on this sequence $(n_k)$, a priori depending on $\mu$.
Investigating partial answers to the Furstenberg conjecture, Lyons asked in \cite{Ly1} whether $S=(s_n)$ could be a {\it universal} such sequence, in the sense that
$\liminf|\widehat\sigma(s_k)|=0$ for every continuous measure $\sigma$ on $\T$ ?

If this holds, the Furstenberg conjecture would be true: if $\mu$ is continuous, then for every $\ell\in\Z$, $\ell\not=0$, $\liminf|\widehat\mu(s_k\ell)|=0$ since $T_\ell\mu$ is still continuous; if in addition $\mu$ is $\times2\times3$-invariant, we get $|\widehat{\mu}(\ell)| = \liminf|\widehat\mu(s_k\ell)|=0$ and $\mu$ is the Lebesgue measure. 
But S. Grivaux and C. Badea \cite{BG} constructed a continuous measure with $\inf_{s\in S}|\widehat\mu(s)|\ge\delta$ for some constant $\delta>0$; in other terms,
\begin{theo}[\cite{BG}] $S$ is a Ka-set. 
\end{theo}
Thanks to an improved version of Drury's result due to Hartman \cite{Har}, we can see that {\it Sidon sets are Ka-sets too}. 

A much more restrictive property of measures has been studied by Bergelson et al.  in the spirit of dynamics (as a spectral property \cite{BJLR}) and by Eisner and Grivaux with an operator-theoretic  point of view (\cite{EG}).
\begin{defi} A subset $E=(n_k)\subset\Z$ is said rigid if there exists a {continuous} probability measure $\mu$ on $\T$ such that $\lim_{k\to\infty}\widehat\mu(n_k)=1$.
\end{defi}
Of course rigid sets are Ka-sets,  whence the question whether  $S$ could be a rigid set ?

\subsection{More on $S$ as a thin set}
\subsubsection{$S$ is $p$-Paley, $1<p<2$}


We mentioned at the beginning of the section that $S$ cannot be $H^1$-Paley. The following is proved in  \cite{GV}: $(L^p, S)$ is a Paley pair. Here is a preliminary result  in this direction, due to Gundy and Varopoulos (cf. \cite{GV}).
\begin{theo}[ \cite{GV}] The set $S$ is  $\Lambda(q)$ for all $q>2$, i.e.
 $||f||_q \le \lambda_q(S)||f||_2$  for any $f\in L^{2}_{S}$.
\end{theo}
The proof uses the {square function}  of a martingale corresponding to a {decreasing} sequence of $\sigma$-subalgebras, and Burkholder's inequalities.
 
\begin{theo}[\cite{GV}]
For any $g\in L^p$ with $1<p<2$,  one has 
$$
\Big(\sum_{n\in S} |\widehat{g}(n)|^2\Big)^{1/2}\leq C_q  \Vert g\Vert_p$$
 where $q>2$ is such that $\frac{1}{p}+\frac{1}{q}=1$, and $C_q<\infty$.  
\end{theo}
\noindent In addition, we proved in \cite{AHM} that  $\lambda_{q}(S)\leq Cq^{3/2}$ and that  $C_q=P(L^p, S)\lesssim q^2$.

\subsubsection{Is $S$ $4/3$-Sidon ?}
We have mentioned at the beginning of the section that $S$ is not  $p$-Sidon for $p<4/3$ and our feeling is that $S$ is $4/3$-Sidon. We can not confirm this yet, but we can prove a partial result
in this direction, by showing that $S$ is $3/2$-Rider. We recall some facts.

If $f=\sum a_k e_k$ is a trigonometric polynomial and $(\varepsilon_k)$ a Rademacher sequence, the randomized polynomial $f_\omega$  is by definition
$$f_{\omega}= \sum \varepsilon_{k}(\omega)  a_k e_k.$$
\begin{defi}\label{prid} A subset $\Lambda$ of $\N$ is called a $p$-Rider set ($1\leq p<2$) if there is a constant $C$ such that, for every $f\in \mathcal{P}_{\Lambda}$, we have 
\begin{displaymath} \Vert \widehat{f}\Vert_p\leq C [f]\end{displaymath} 
where $[f]$ denotes the Pisier norm of $f$ defined by
$$[f]=\mathbb{E}(\Vert f_{\omega}\Vert_\infty).$$ 
\end{defi}

The $1$-Riderness  coincides with the  $1$-Sidonicity, i.e. Sidonicity (cf. \cite{Ri}).  For $1<p<2$,  the $p$-Sidonicity  (namely $\Vert \widehat{f}\Vert_p\leq C \Vert f\Vert_\infty$)
implies the $p$-Riderness,  the converse being  open. But it  is more flexible to work with  $p$-Riderness than  with $p$-Sidonicity. See \cite{AHM} for more on this.

In order to study the $p$-Riderness of $S$, we invoke two theorems of  Rodr\'i\\guez-Piazza (\cite{Ro1}, Lema 2.4 p. 89 and Teorema 2.3 p. 85-86), which were used in \cite{LQR}.
For stating these theorems, we need the following notations.
For an arbitrary  subset $\Lambda$ of $\N$ and $n=1,2,\ldots$, let us set once and for all
\begin{displaymath}\Lambda_n=\Lambda\cap [1,n] \hbox{\ and}\  \Lambda_{I_n}=\Lambda\cap I_n \hbox{\quad where}\  I_n=[2^n, 2^{n+1}[\subset \N. \end{displaymath}
For a finite subset $A$ of $\N$, we write 
$$\psi_A=\big\Vert \sum_{k\in A} e_k\big\Vert_{\psi_2}.$$
Recall that $\psi_2$ designates the gaussian Orlicz function: $\psi_{2}(x)=e^{x^2}-1$ and that, for a function $f$, it holds  \cite[p.44, vol.1]{LQ}
\begin{equation}\label{psi2}
\Vert f\Vert_{\psi_2}\asymp \sup_{q\geq 2} \frac{\Vert f\Vert_q}{\sqrt{q}}.
\end{equation}

The following result gives a lower bound for the size of the largest quasi-independent  subset in a given finite set.

\begin{theo} [\cite{Ro2}, p.89]\label{fcs} Let $A\subset \mathbb{N}$ be a finite set. We can find a quasi-independent  subset $E$ of $A$ with cardinality 
$$|E|\geq \delta\Big(\frac{|A|}{\psi_{A}}\Big)^2$$
where $\delta$  is a positive constant.
  \end{theo}
  
 The following is a necessary and sufficient condition for a set $E\subset \mathbb{N}$ to be $p$-Rider, a condition involving the size of 
  the largest quasi-independent  subset in an arbitrary finite subset of $E$.

\begin{theo}  [\cite{Ro2}, p.85] \label{betis} Let $ 1\le p<2$.  A subset $E\subset \mathbb{N}$
  is $p$-Rider if and only if for every finite subset $A$ of $E$, there exists a quasi-independent subset $B$ of $A$ such that 
$$  |B|\geq \delta |A|^\varepsilon$$
 where $0<\varepsilon=\frac{2}{p}-1\leq 1$  and where  $\delta$ is a positive constant.
\end{theo}

The case $p=1$ is due to Pisier. Observe that $\frac{2}{p}-1=1/2$ when $p=4/3$. \\
We begin with a general fact.
\begin{theo} \label{gene} Let $\Lambda\subset \N$. Suppose  there exist constants $C>0$ and $\alpha\geq 1/2$
such that 
$$\forall q>2,\ \forall f\in \mathcal{P}_{\Lambda},  \ \   \Vert f\Vert_q \leq C q^\alpha \Vert f\Vert_2.$$
Then 
$\Lambda$ is $p$-Rider for $p=\frac{4\alpha}{2\alpha+1}\cdot$
 \end{theo}
 \begin{proof}
 Let $A$ be a finite subset of $\Lambda$ with $|A|=n$. We claim that 
\begin{equation}\label{dim}\Vert \sum_{j\in A} e_j\Vert_q \lesssim \min( q^{\alpha}\sqrt{n},\  n)\hbox{\ for all}\ q\geq 2.\end{equation}
Indeed, firstly, we have obviously  $\Vert \sum_{j\in A} e_j\Vert_q\leq \Vert \sum_{j\in A} e_j\Vert_\infty= n$; secondly the assumption allows us  to get 
$$\Vert \sum_{j\in A} e_j\Vert_q \lesssim  q^{\alpha}\Vert \sum_{j\in A} e_j\Vert_2= q^{\alpha} \sqrt{n}.$$ 
From (\ref{dim}) we easily get  (separating the cases $q\leq n^{\frac{1}{2\alpha}}$,\  $q\geq n^{\frac{1}{2\alpha}}$) that 
$$\psi_{A}\asymp \sup_{q\geq 2}\frac{1}{\sqrt{q}}\Vert \sum_{j\in A} e_j\Vert_q \lesssim n^{1-\frac{1}{4\alpha}},$$
using the formula (\ref{psi2}).
Now,  Theorem \ref{fcs}  provides us with a quasi-inde\\pendent subset $E$ of $A$ of cardinality 
$$
   |E|\gtrsim \left(\frac{n}{\psi_A}\right)^2\gtrsim n^{\frac{1}{2\alpha}}.
   $$
Finally, adjust the exponent  $p$ so as to have
$$\frac{1}{2\alpha}=\frac{2}{p}-1,$$ that is $p=\frac{4\alpha}{2\alpha+1}$. Then Theorem \ref{betis} allows us to conclude.
\end{proof}

If $\alpha=1/2$,  we have $\frac{4\alpha}{2\alpha+1}=1$, and $\Lambda$ is a Sidon set.
  
\medskip

The Furstenberg set $S$ is $\Lambda(q)$ for all $q<\infty$  with specified constants \cite[Theorem 4.12]{AHM}). This result allows us  to prove 
 a partial result concerning the  $p$-Riderness of $S$. It is partial  because we conjecture that $S$ is $4/3$-Rider (perhaps even $4/3$-Sidon), but we only prove its $3/2$-Riderness.
 More generally, we conjecture that  $S(q_1,\ldots, q_s)$   (where the $q_j$'s are multiplicatively independent integers) is $\frac{2s}{s+1}$-Rider (or even $\frac{2s}{s+1}$-Sidon) and this would be optimal. We prove the following partial
  result for $S(q_1,\ldots, q_s)$.

\begin{theo} \label{decept} The Furstenberg  set $S$ is $3/2$-Rider.
 More generally, let $s\ge 2$, we have
 \begin{enumerate}
 \item  The set $S(q_1,\ldots, q_s)$ is $\frac{2s-1}{s}$-Rider. 
 \item
 The set  $S(q_1,\ldots, q_s)$ is not $p$-Rider for $p<\frac{2s}{s+1}$.
 \end{enumerate}
 \end{theo}
\begin{proof}
We appeal to  \cite[Theorem 4.12]{AHM} which tells that the assumption of Theorem \ref{gene} holds for $\Lambda=S(q_1,\ldots, q_s)$ with 
$\alpha=s-1/2$,  giving
 $$\frac{4\alpha}{2\alpha+1}=\frac{2s-1}{s}\cdot$$
  The multiplicative character of $S(q_1,\ldots, q_s)$ lurks in the value of that exponent $\alpha$. The second assertion comes from the mesh condition for $p$-Rider sets \cite[Proposition 3.2]{AHM}: if $\Lambda$ is $p$-Rider, we must have 
 $|\Lambda_N|\lesssim (\log N)^{p/(2-p)}$. 
 But we know  \cite[page 9]{AHM} that $|\Lambda_N|\asymp (\log N)^s$. So that 
 $s\leq p/(2-p)$, or again  $p\geq\frac{2s}{s+1}$.
\end{proof}

\subsubsection{$S$ and the Bohr topology}
Recall that if  $G$ is a locally compact abelian group with dual $\Gamma$, the  Bohr  compactification $\beta\Gamma$ of  $\Gamma$ is the dual group of $G_d$, the group $G$ equipped with the discrete topology. The group $\beta\Gamma$ is the set of all characters (continuous or not) on $G$, it is compact and contains $\Gamma$ as a dense subgroup. We describe this topology when $G=\T$ and  $\Gamma=\Z$.
\begin{defi} The Bohr topology on  $\Z$ is  the  group topology with the following basis of neighbourhoods of  zero: 
$$V(x_1,\dots,x_k,\eps)=\{n\in\Z; \ |e(nx_j)-1|<\eps\,  \hbox{\ for}\ 1\le j\le k\}$$
with $\varepsilon>0$ and $ x_j\in\R$,
called \textnormal {Bohr neighbourhoods} (of $0$ in $\Z$).
\end{defi}
The Bohr topology on  $\Z$ is the coarsest  topology for which all Fourier transforms of discrete measures on $\T$ are continuous.

The following assertion follows from the pigeonhole principle and the simultaneous diophantine approximation.
\begin{prop} A Bohr neighbourhood  has positive upper density.
\end{prop}

%

We are concerned with  $\beta\Z$ 
and with some subsets of integers, dense in $\beta\Z$ (i.e. Bohr-dense) or not.  A first class of examples is the following.
\begin{prop} A Hartman-distributed set is Bohr-dense.
\end{prop}
It is a direct consequence of the definition.
This property of Bohr-density appears in the dynamical classification of (rather big) sets of integers intensively studied by Bergelson, Bourgain, Ruzsa and many others (\cite{BL,Bo1,Ru}). We extract the main implications we need:
$$
 {\rm Hartman-distributed} \Longrightarrow {\rm  Uniformly \ recurrent}   \Longrightarrow {\rm Bohr-dense},
 $$
with the definitions to come.
\begin{defi} 
A set $E\subset\Z$ is said to be recurrent if, for every dynamical system $(X,\mathcal{A},\mu,T)$ and every subset $A\in \mathcal{A}$ with positive measure, there exists $h\in E$ such that $\mu(A\cap T^{-h}A)>0$.
\end{defi}
\begin{defi} A set $E\subset\Z$ is said to be uniformly recurrent if every translate of $E$ is recurrent.
\end{defi}
The question whether the last implication above is reversible remains open.

\smallskip

Back to the Furstenberg set, we show that
\begin{theo}\label{BC} The set $E:=S(p_1,\dots, p_r)$, in particular $S$, is Bohr-closed, thus cannot be a recurrent set. 
\end{theo}
\begin{proof} We write  $p^\alpha$ for $p_{1}^{\alpha_1}\cdots p_{r}^{\alpha_r}$, and we write  $\beta\geq \alpha$ if $\beta_j\geq \alpha_j$ for all  $j$; and  $\beta>\alpha$ if $\beta_j > \alpha_j$ for all $j$. We now show that $\Z\backslash E$ is open for the  Bohr topology.
Indeed, let $m\notin E$. We distinguish two cases.  \\
\noindent $\bullet$  $m=0$. Then $V=N\Z$, where $N$ has a prime factor  $p>p_r$, is a neighbourhood  of  $0$ disjoint from $E$.\\
\noindent $\bullet$  $m\neq 0$. 
One  writes 
$$m=p^\alpha n=:s_0\times n$$ 
with $s_0\in E$,\  $n\not=0,1$ and $n\wedge  p_{1} p_{2}. ..p_r=1$. 
 One can find $s=p^\beta\in E$ with $\beta>0$ and the $\beta_j$'s large enough so as to have $s>n-1$, implying that   $n\not\equiv1$ mod $s$. 
We now  claim that the neighbourhood  of $m$,
$$V:=V(m)=m+sp^\alpha \Z,$$ 
satisfies $V\cap E=\emptyset$. Indeed, a  relation
$$p^\gamma=m+sp^\alpha k= p^\alpha(n+sk),$$
clearly implies $\gamma\geq \alpha$. If $\gamma=\alpha$, then $n+sk=1$ and $n\equiv1$ mod $s$, contradicting the choice of  $s$. Therefore, we have for example  $\gamma_1>\alpha_1$. 
 After simplification by $p_{1}^{\alpha_1}$, we get
$p_1|(n+sk)$ which is again  impossible: since $p_1|s$, we have 
 $(n+sk)\wedge p_1=n\wedge p_1=1$. 
\end{proof}

 A second dynamical classification due to Hartman concerns small sets of integers and it seems to escape this Bohr property: 
 $$ E \ {\rm  Sidon\ or} \ E\ {\rm  rigid}  \Longrightarrow E\ {\rm Ka-set} \Longrightarrow W(E) \  {\rm uncountable}.
 $$ 
 However, Katznelson \cite{Kat} constructed a Bohr-dense Ka-set and Griesmer \cite{Gr} constructed a both rigid and Bohr-dense set. The existence of a Sidon set dense in $\beta\Z$ remains open.

\section{A random version $T$ of $S$}
Recall that we define a random version $T$ of $S$ by
$$T=\{k\in \N: \xi_k=1\}$$
where $(\xi_k)$ is a sequence of independent $0$-$1$-valued independent random variables with $\mathbb{E}(\xi_k)=\delta_k:=\frac{ \log k}{k}$.
\subsection{First comparative results}
Here are some results about the random version $T =(t_n)$.  
First, $T$ shares with $S$ some {sparseness properties.}

\begin{theo}[\cite{AHM}]\label{tijdeman} Almost surely, the difference $t_{n+1}-t_n$ satisfies\textnormal{:}\\ 
\indent {\rm 1)}  $\displaystyle\limsup_{n\to \infty} \frac{t_{n+1}-t_n}{(t_n/\log t_n)\,\log\log t_n}\leq 2;$\\
\indent {\rm 2)} $\displaystyle\liminf_{n\to \infty} \frac{t_{n+1}-t_n}{t_n/(\log t_n)^{3+\delta}}\geq 1 \hbox{\quad for all}\ \delta>0.$\\
In particular, the set $T$ satisfies: {$t_{n+1}-t_n\to \infty,\ t_{n+1}/t_n\to 1 $\  a.s.}
\end{theo}
Here is a {dynamical property} of $T$, shared by more general random sets $R$ studied by Bourgain. Let  
$$m_N:=\sum_{k=1}^N \delta_k.
$$ 

\begin{theo}[\cite{Bo2}]  If $\delta_k \downarrow$ and $m_N/\log N\to \infty$, the corresponding random set $R$  is almost surely {Hartman distributed.} In particular,  {that is the case for   $T$.}
\end{theo} 
This property of $R$ is in strong contrast with the case of $S$.
 \medskip
 
Now comes a {harmonic analysis property} of $T$, which is better than  that of  $S$.
\begin{theo}[\cite{LQR2}, \cite{AHM}] 
Almost surely, the set $T$ is  $p$-Rider for $p>4/3$ and  not $p$-Rider for $p<4/3$.  
\end{theo} 
The case $p<4/3$ is trivial {(mesh condition).} We do not know if $T$ is $4/3$-Rider but we tend towards this result in the next item.

\subsection{A large subset $T'$ of $T$ is $4/3$-Rider} 
In the sequel, we prove the $4/3$-Riderness and $\Lambda(q)$ property of a  large subset $T'$ of $T$.  The detailed proof completes a highly sketchy proof in \cite[Lemma 3.3]{LQR2}. We do not know whether  $T$ itself  possesses these properties of $T'$.

 We begin with some notations. Let $E\subset \mathbb{N}$.
 For $E'\subset E$, the upper  density of $E'$ inside $E$ is defined by
$$\overline{d} (E',E)=\limsup_{N\to \infty} \frac{|E'\cap [1,N]|}{|E\cap [1,N]|}.$$
We define similarly the lower density $\underline{d}(E',E)$ with $\liminf$.  Clearly
$$0\leq \underline{d}(E',E)\leq \overline{d}(E',E)\leq 1.$$
Clearly, $\underline{d}(T',T)>0$ indicates that $T'$ occupies a large portion of $T$.

\begin{theo} \label{tiers} Almost surely, the random set $T$ contains a subset $T'$ of \textnormal{positive lower density} such that
\begin{enumerate} 
 \item $T'$ is a  $\frac{4}{3}$-Rider set;
 \item
 $T'$ is  a $\Lambda(q$)-set for all $q<\infty$.
 \end{enumerate}
\end{theo}

 The key point for proving  Theorem \ref{tiers} is the our main technical result below, which is stated with notations of Section 3.2.2
 and the long proof of which is outlined at the end of this subsection.
 
\begin{theo} \label{quarte}   Almost surely, the following two estimates hold:
\begin{enumerate}
\item  $|T_{I_n}|\asymp n$.
\item  $\psi_{T_{I_n}}\leq C|T_{I_n}|^{1/2}$ where $C=C(\omega)$ does not depend on $n$. 
\end{enumerate}
\end{theo}

\subsubsection{Proof of Theorem \ref{tiers}.}
Taking Theorem \ref{quarte} for granted, we can prove Theorem \ref{tiers} by a combinatorial argument of \cite{LQR} which we reproduce here. First, by Theorem \ref{fcs}
and Theorem \ref{quarte} (2) ,  we can find a quasi-independent  subset $E_n$ of $T_{I_n}$ with cardinality 
$$|E_n|\geq \delta\Big(\frac{|T_{I_n}|}{\psi_{T_{I_n}}}\Big)^2\geq \delta' |T_{I_n}|.$$
Let now $$T':=\bigcup_{n\geq 1} E_n.$$
Notice that  $T'_{I_n}=E_n$ is a quasi-independent subset  of   size proportional to $|T_{I_n}|$.  We are going to see that something similar subsists for \textit{arbitrary finite} subsets of $T'$,  with some loss.

We first note that $T'$ has positive lower density in $T$.  Indeed, if $N$  is given and $2^{n}\leq N<2^{n+1}$, we know that almost surely  $|T_N|\lesssim (\log N)^2$ while 
$$|T'_N|\geq \sum_{k=1}^{n-1} |E_k|\geq \delta' \sum_{k=1}^{n-1} |T_{I_k}|\gtrsim \sum_{k=1}^{n-1} k\gtrsim n^2\gtrsim (\log N)^2\gtrsim |T_N|,$$
where for the first $\gtrsim$ we have used Theorem \ref{quarte} (1).

 Let now  $A\subset T'$ be an arbitrary finite set. We claim that $A$ contains a quasi-independent subset $E$ of size $\gtrsim |A|^{1/2}$. To this end, we put 
$$J=\{n: A\cap E_n\neq \emptyset\}=\{n_1<n_2<\cdots< n_h\}$$
so that
$$
 A=\bigcup_{n\in J} (A\cap E_n).
$$
We look for a big quasi-independent subset $B$ in $A$ such that $|B\gtrsim |A|^{1/2}$  by distinguishing  two cases: \\

\noindent {\it Case I. $|A\cap  E_n|\geq |A|^{1/2}$ for some $n\in J$}. Then, take  $B:=A\cap  E_n\subset  E_n$.\\

\noindent {\it Case II.   $|A\cap  E_n|\leq |A|^{1/2}$ for all $n\in J$.} Then, $h\geq |A|^{1/2}$. Pick a point $\mu_j$ in $A\cap I_{n_{2j+1}}$ for each $j\leq K:= \big[(h-1)/2\big]$ and observe that $\mu_{j+1}/\mu_j\geq 2$ so that $B:=\{\mu_1,\mu_2,\ldots, \mu_K\}$ is a quasi-independent set with cardinality $\geq \delta'' h\geq  \delta'' |A|^{1/2}$. 
 This proves our claim. 
 
 As the assumptions of Theorem \ref{betis} are satisfied for $T'$,  with $p$ such that $(2/p)-1=1/2$, i.e. $p=4/3$, the first part of Theorem \ref{tiers} is thus proved.
 
Let us now turn to the second part of Theorem \ref{tiers}.  Assume $f\in L^{2}_{T'}$ and let
 $f_n=\sum_{k\in  E_n} \widehat{f}(k)e_k$.   The Littlewood-Paley theorem \cite[Chapter 8]{MS}  and the convexity of the $L^q$-norm  for $q\geq 2$ imply
 $$\Vert f\Vert_q\leq C_q\Vert  \big(\sum_{n} |f_n|^2\big)^{1/2}\Vert _q\leq C_q\big(\sum_{n} \Vert f_n\Vert_{q}^{2}\big)^{1/2}
 $$
 with $C_q\lesssim q^{3/2}$ \cite{LQR}.
 But $ E_n$ is quasi-independent, hence Sidon with a Sidon constant $\leq 8$; consequently, $E_n$ is a $\Lambda(q)$-set with $\lambda_{q}( E_n)\leq C\sqrt{q}$, so that $ \Vert f_n\Vert_{q}\leq C \sqrt{q}  \Vert f_n\Vert_{2}$
 and we get 
 $$\Vert f\Vert_q\lesssim C_q \sqrt{q}\big(\sum_{n} \Vert f_n\Vert_{2}^{2}\big)^{1/2}= C_q \sqrt{q} \Vert f\Vert_2.$$ 
  This means that $T'$ is a $\Lambda(q)$-set, with  $\lambda_{q}(T')\lesssim q^2$.
   $\Box$
   \medskip
    \subsubsection{Proof of $\psi_{T_{I_n}}\leq C|T_{I_n}|^{1/2}$}
 We now prove the second part of Theorem \ref{quarte} (the first one is easy). In what follows, $C_0,C_1,\ldots $ will denote constants. We begin with a simple interpolation lemma--estimating
 $\|f\|_{\psi_2}$ by $\|f\|_\infty$ and $\|f\|_2$. 
\begin{lem}\label{interpol} Let $f\in L^\infty$. Then
\begin{equation}\label{interpola} \Vert f\Vert_{\psi_2} \leq \frac{\Vert f\Vert_\infty}{\sqrt{\log (1+\Vert f\Vert_{\infty}^{2}/\Vert f\Vert_{2}^{2})}}.\end{equation}
In particular, if $f=\sum_{k\in A} c_k e_k$ with $A$ finite and $c_k$ scalars, we have
\begin{equation}\label{majorant} \Vert f\Vert_{\psi_2}\leq \frac{\sum |c_k|}{\sqrt{\log \big(1+\frac{(\sum |c_k|)^{2}}{\sum |c_k|^2}}\big)}.\end{equation} 
\end{lem}
\begin{proof} We can assume $\Vert f\Vert_\infty=1$. Let $\lambda= \Vert f\Vert_{\psi_2}$ and $\psi_{1}(x)=e^{x}-1$. By the definition of $\lambda$  and the convexity of  $\psi_1$ together with the facts  $|f|\leq 1$ and $\psi_{1}(0)=0$, we have
$$1=\int \psi_{2}(|f|/\lambda)=\int \psi_{1}(|f|^2/\lambda^2)\leq  \psi_{1}(1/\lambda^2)\int |f|^2.$$
Inverting this relation gives (\ref{interpola}), since $\psi_{1}^{-1}(y)=\log (1+y)$.

  Next, observe that 
  $x\mapsto x/\log(1+x)$ increases on $\R^+$.  Then 
   (\ref{interpola}) implies    (\ref{majorant}) through 
 the relations 
 $\Vert f\Vert_\infty\leq \sum_{k\in A}|c_k|,\  \Vert f\Vert_{2}^{2}=\sum_{k\in A} |c_k|^2.$  
 \end{proof}
 
 The following specializations of Lemma \ref{interpol} will be used.
 \begin{lem}\label{weakmo} Let $f=\sum_{k\in I_n} a_k b_ke_k$ with $a_k, b_k$ scalars and $\sum |a_k|^2\leq 1$.
 Then
\begin{equation}\label{part}\Vert f\Vert_{\psi_2}\leq  \frac{(\sum |b_k|^2)^{1/2}}{\sqrt{\log (1+(\sum |b_k|^2)^{1/2}})}.
\end{equation} 
Moreover, with  $\delta_k=\frac{k}{\log k}$, 
\begin{equation}\label{inpart}\Vert \sum_{I_n} \delta_k e_k\Vert_{\psi_2}\leq C_0 \sqrt{n}.
\end{equation}
\end{lem}
\begin{proof}
Apply (\ref{majorant}) with  $c_k=a_kb_k$  
and Cauchy-Schwarz inequality, remembering that \ $\sum |a_k|^2\leq 1$. The maximum of $(\sum |c_k|)^2/\sum|c_k|^2$ is obtained when $\sum |a_k|^2=1$ and $a_k$ is proportional to  $|b_k|$,  
that is  $a_k={|b_k|}/\sqrt{\sum |b_k|^2}$.  Then $\sum |c_k|=\sqrt{\sum |b_k|^2}$, and this gives (\ref{part}).

 Next, for  (\ref{inpart}), we observe that
$\sum_{k\in I_n} \delta_k\asymp 2^n\times (n/2^n)\asymp n$  and similarly $\sum_{k\in I_n} \delta_{k}^{2} \asymp  2^n\times (n^{2}/4^n)\asymp  (n^{2}/2^n)$, so that ${(\sum \delta_k)^2}/{\sum\delta_{k}^{2}}\asymp 2^n$. 
Now, (\ref{majorant}) gives
$$\Vert \sum_{I_n} \delta_k e_k\Vert_{\psi_2}\lesssim  \frac{n}{\sqrt{\log 2^n}}\asymp \sqrt{n}. $$
\end{proof}
We will apply the above lemma to the random function $f=  \sum_{k\in I_n} (\xi_k-\delta_k) e_k$, considered as function in the Banach space $L^{\psi_2}$.
 For this, we first need a simple symmetrization lemma. 
\begin{lem}\label{cartesym} Let  $\mathcal{X}$ be a Banach space, and $X$ be a $\mathcal{X}$-valued and integrable random variable with symmetrization $\tilde{X}=X-X'$ where $X'$ is an independent copy of $X$. Let $t>2\,\mathbb{E}(\Vert X\Vert)$. Then
\begin{equation}\label{toto}  \mathbb{P}(\Vert X\Vert>2t)\leq 2  \mathbb{P}(\Vert \widetilde X\Vert>t).\end{equation}
\end{lem}
\begin{proof} Firstly, Markov's inequality gives  
$$ \mathbb{P}(\Vert X\Vert>t)\leq \frac{\mathbb{E}(\Vert X\Vert)}{t}\leq \frac{1}{2}.$$
Secondly,  $ \Vert X\Vert>2t$ and   $\Vert X'\Vert\leq t$ imply $\Vert \widetilde X\Vert>t$, so that
\begin{eqnarray*}
\mathbb{P} ( \Vert \widetilde X\Vert>t)
&\geq&\mathbb{P}(\Vert X'\Vert\le t, \Vert X\Vert>2t) \\
&=&\mathbb{P}( \Vert X'\Vert\leq t) \mathbb{P}(\Vert X\Vert>2t)
\geq  \frac{1}{2}\mathbb{P}(\Vert X\Vert>2t).
\end{eqnarray*}
Then (\ref{toto}) follows.\end{proof}

Consider  $(X_k)_{k\in I_n}$ with $X_k=\xi_k -\delta_k$ and their symmetrizations $(\widetilde{X}_k)_{k\in I_n}$. 
Let now
$$Z_n=\Vert \sum_{k\in I_n} X_k e_k\Vert_{\psi_2}, \quad \widehat{Z}_n=\Vert \sum_{k\in I_n}\widetilde{X}_k e_k\Vert_{\psi_2}.$$ 
 Thanks to (\ref{inpart}), 
 it holds
\begin{equation} \label{start} 
\psi_{T_{I_n}}:=\Vert \sum_{k\in I_n} \xi_k e_k\Vert_{\psi_2}\leq Z_n+\Vert \sum_{I_n} \delta_k e_k\Vert_{\psi_2}\leq Z_n+C_0\sqrt{n}.
\end{equation} 

  In order to  majorize $\widehat{Z}_n$ and then $Z_n$, we need an  inequality, proved independently by several authors,  
    often referred to as Talagrand's deviation inequality for Lipschitz functions (see   e.g.  \cite{BLM},  \cite[Corollary 4 p.~75]{JoSc},  \cite[Theorem 3]{Tal}). 
   Here is a version borrowed from \cite{JoSc}.
  \begin{theo}[\cite{JoSc}, p.75] \label{gide}
  Let $f:\ell_{2}^{n} \to \R$ be a convex function with Lipschitz constant $\lambda$. Let 
  \begin{displaymath} Z=f(\varepsilon_1,\ldots, \varepsilon_n) \end{displaymath} 
  where $(\varepsilon_j)_{1\leq j\leq n}$ is a Rademacher sequence. Then, it holds
  $$\mathbb{P}\big(|Z-\mathbb{E}(Z)|>t\big)\leq a\exp(-b\frac{t^2}{\lambda^2}) \hbox{\quad for all}\  t\geq 1$$
  where $a,b$ are positive absolute constants.
   \end{theo}
  As a corollary, we get
    \begin{theo}\label{mofa} Let $A\subset \N$ be a finite set. Let $(\varepsilon_j)_{j\in A} $ be a Rademacher sequence,  and $v=(v_j)_{j\in A}$ be vectors in a Banach space $\mathcal{X}$, with \textnormal{weak moment} $\sigma$ defined by
\begin{equation}\label{haak}\sigma=\sigma(v):=\sup_{\varphi\in B_{\mathcal{X}^\ast}}\big(\sum_{j\in A} |\varphi(v_j)|^{2}\big)^{1/2}=\sup_{\sum |a_j|^2\leq 1} \Vert \sum_{j\in A} a_j\,v_j\Vert\end{equation}
where $B_{\mathcal{X}^\ast}$ is the unit ball of $\mathcal{X}^\ast$.
Set $Z=\Vert \sum_{j\in A} \varepsilon_j v_j\Vert$.  Then, for all $t>0$, the following \textnormal{two-sided} inequality holds with absolute constants $a,b>0$:
\begin{equation}\label{salaud}  \mathbb{P}\big(|Z-  \mathbb{E}(Z)|>t\big) \leq a \exp\big(-b \frac{t^2}{\sigma^2}\big).\end{equation} 
\end{theo}
\begin{proof} We apply Theorem \ref{gide} to the convex function $f(x)=\Vert \sum_{j\in A} x_j v_j\Vert$.
It suffices to note that $f$ has Lipschitz constant $\lambda$ exactly equal to $\sigma$, which is elementary via Hahn-Banach's theorem.
\end{proof}
 

To exploit this theorem, it is convenient to note first  the following. 
\begin{lem}\label{khik} For $Z_n=\Vert \sum_{k\in I_n} X_k e_k\Vert_{\psi_2}$, we have 
$$\mathbb{E}(Z_n)\leq \ \mathbb{E}(\widehat{Z}_n)\leq C_1\sqrt{n}.$$
\end{lem}
\begin{proof} The left inequality is clear since the $X_k$'s are centered \cite[Theorem 2.6]{Hof}. Next, the symmetry of the $\widetilde{X}_k$, the $L^2$-$L^{\psi_2}$ Khintchine inequalities and Fubini's theorem   imply 
$$\mathbb{E}(\widehat{Z}_n)\lesssim \big(\sum_{k\in I_n} V(X_k)\big)^{1/2} = (\sum_{k\in I_n} \delta_k(1-\delta_k)\big)^{1/2} \leq C_1\sqrt{n}.$$
\end{proof}  
We next claim that
\begin{equation}\label{titi}  t>2C_1\sqrt{n}\Longrightarrow \mathbb{P}(Z_n>4t)\leq \mathbb{P}\Big(\widehat{Z}_n-\mathbb{E}(\widehat{Z}_n)>t\big).\end{equation}
Indeed, since $t>2C_1\sqrt{n}\geq 2 \mathbb{E}(Z_n)$, the relation (\ref{toto}) gives us 
$$\mathbb{P}(Z_n>4t)\leq 2 \mathbb{P}(\widehat{Z}_n>2t)\leq  2\mathbb{P}(\widehat{Z}_n-\mathbb{E}(\widehat{Z}_n)>t)$$
since $\mathbb{E}(\widehat{Z}_n)\leq C_1\sqrt{n}\leq t$.
 
 Hence,  we are led to find upper bounds on the RHS of (\ref{titi}).
 We claim that 
\begin{equation}\label{oneget} \mathbb{P}\big(\widehat{Z}_n-\mathbb{E}(\widehat{Z}_n)>t\big)\leq a \int_{\Omega} \exp\big(-b \frac{t^2}{\sigma_{\omega}^2}\big) d\omega\end{equation}
where $\sigma_{\omega}$ denotes the  weak moment of the vectors $v_k= \widetilde{X}_{k}(\omega) e_k,\  k\in I_n$, in the Banach space $\mathcal{X}=L^{\psi_2}$.
Indeed, 
it suffices to apply  Theorem \ref{mofa} to the variables $\varepsilon_k$ and the vectors $X_kv_k$, and  a symmetrization argument.

In the following, we estimate $\sigma_\omega$. For simplicity,  we will abbreviate $\sum_{k\in I_n}$ to $\sum$ and $\prod_{k\in I_n}$ to $\prod$. 
First, remark that Lemma \ref{weakmo} implies 
\begin{equation} \label{sigma} \sigma_\omega \leq  \frac{(\sum |\widetilde{X}_k|^2)^{1/2}}{\sqrt{\log (1+(\sum |\widetilde{X}_k|^2)^{1/2}})}.\end{equation}
 Second, we are going to prove that 
 $$
    W_n: =\sum_{k\in I_n} |\widetilde{X}_k|^2=\sum |\xi_k-\xi_k'|^2 = O(n)
 \quad a.s. 
 $$ 
 by showing  
 \begin{equation}\label{tau} 
    \mathbb{P}(W_n \geq 5C_1 n) \leq \exp(-C_1 n).
    \end{equation}
 Indeed, since  $|\widetilde{X}_k|=|\xi_k-\tilde{\xi}_k|\leq 1$ takes only values $0$ and $1$ with $P(|\widetilde{X}_k|=1)=2\delta_k(1-\delta_k)$,
 Markov's inequality implies  
$$
  \mathbb{P}(W_n \geq 5C_1 n)
    \leq 
  e^{-5C_{1} n} \prod_{}  \mathbb{E}(e^{ |\widetilde{X}_k|^2})
  $$
  But
  $$
 \prod \mathbb{E}(e^{ |\widetilde{X}_k|^2}) \le
 \prod_{} \mathbb{E}(1+2 |\widetilde{X}_k|^2)
 \leq   \prod_{}(1+4\delta_k)\leq e^{\sum_{} 4\delta_k}\le
 e^{4C_{1}n}. 
$$
We have thus proved \eqref{tau}.

 If $W_n(\omega) < 5C_1 n$,  we have $\sigma_\omega \leq C_2 \sqrt{{n}/{\log n}}$ according to (\ref{sigma}).  This,  together with (\ref{tau}), allows us to write  (\ref{oneget}) under the form 
 \begin{equation}\label{form}  \mathbb{P}\big(\widehat{Z}_n- \mathbb{E}(\widehat{Z}_n)>t\big)\leq a\int_{W_n <5C_1 n}\!\exp\big(-b \frac{t^2}{\sigma_{\omega}^2}\big) d\omega + a \mathbb{P}(W_n\ge 5C_1 n)
 $$
$$\leq a \exp(-b \frac{t^2}{C_{2}^{2}}\frac{\log n}{n})+ a \exp (-C_1 n). 
 \end{equation} 
 Taking $t=t_n=C_3 \sqrt{n}$ with large $C_3$ 
 so that ${b\,C_{3}^{2}}/{C_{2}^{2}}\geq 2$, we get in view of (\ref{titi}) and (\ref{form}):
  \begin{equation*}\label{forme}  
  \sum \mathbb{P}(Z_n>4t_n)\leq  a\sum \Big( n^{-2}+ e^{-C_1n}\Big)<\infty.
  \end{equation*}  
 Now, by Borel-Cantelli's lemma, there exists almost surely an integer $n_0=n_{0}(\omega)$ such that, for all $n\geq n_0$, it holds 
  $Z_n\leq 4t_n= 4C_3\sqrt{n} $.  Hence, since we have as well $|T_{I_n}|\asymp n$ almost surely, we get 
  $$\psi_{T_{I_n}} \leq Z_n+C_0\sqrt{n} \leq (4C_3+C_0)\sqrt{n}\leq C_4|T_{I_n}|^{1/2},$$
 This ends  the proof of Theorem \ref{quarte}.  
  \hfill$\square$
  \medskip

  
\subsection{Sharpness of Bourgain's random condition}
 Let $R$ be the random set of integers associated to a Bernoulli sequence $(\xi_k)$ with $\mathbb{E}(\xi_k)=\delta_k\downarrow$ and let $m_N=\sum_{k=1}^N \delta_k$. We saw in \cite[Theorem 5.3]{AHM} that  (Bourgain's result) if $m_N/\log N\to \infty$, $R$ is almost surely H-distributed. We mention in passing that a result of \cite{Hu} is wrong: the same conclusion is claimed when $m_N\to \infty$. But this \textit{cannot work:}   a (correct) result of Kahane and Katznelson \cite[Theorem 1 p.~364]{KK1}  shows that if $\delta_k=1/k$,  the corresponding $R$ is almost surely Sidon, hence not Hartman-distributed. We are going to show here a kind of converse which on the one hand shows that Bourgain's result is rather sharp, and on the other hand  improves on Theorem 1 in \cite{KK1}, which was proved by a method involving the theory of multiplicative chaos.

  \begin{theo}\label{mieux que kk} Assume that $\delta_k$ is decreasing and $m_N=O(\log N)$. Then $R$ is almost surely Sidon, hence not Hartman-distributed.
   \end{theo}
\begin{proof} We proceed in two steps.\\
 {\bf Step 1.} {\it We assume that 
 $m_N\leq c\log N$  for large $N$ \hbox{\ with \ some}\ $c\leq 1/(24 e)$}. \\ 
  We first recall a lemma from \cite{LQR} already used in \cite{AHM}.
\begin{lem} \label{trio} Let $n\geq 2$ and $A\ge 1$ be positive integers.  Set
$$\Omega_{n}(A)= \{\omega : R(\omega)\cap [A,\infty[ \hbox{\ contains at least a relation of length}\  n \}.$$
 Then
 $$\mathbb{P}(\Omega_{n}(A))\leq \frac{B^n}{n^n} \sum_{j>A} \delta_{j}^2 m_{j}^{n-2}, \hbox{\ with}\ B=4e.$$
\end{lem}

Recall that $R_k=R\cap [1,k],\  k=1,2,\ldots $.  We will prove 
\begin{lem}\label{nouv} Let $A_n=e^{n}$.  Then\\
\indent {\rm (1)} $\sum_{n\geq 1} \mathbb{P}(\Omega_{n}(A_n))<1.$\\
\indent {\rm (2)} Almost surely,  $|R_{A_n}|\leq n$ for all integers $n$ large enough. 
\end{lem}
Indeed, since $n\delta_n\leq \sum_{k=1}^n \delta_k\leq c\log n$, we get from Lemma \ref{trio} that 
$$\mathbb{P}(\Omega_{n}(A_n))\leq \frac{B^n c^n}{n^n} \sum_{j>A_n}\frac{(\log j)^n }{j^2}.
$$
Let $f_n(t) ={(\log t)^n}/{t^2}$. This function decreases on $[A_n, \infty[$ since
$$f'_{n}(t)=\frac{(\log t)^{n-1}}{t^3} \big(n-2\log t\big)\leq 0$$
  when $t>A_n=e^n$, and then
$$\sum_{j>A_n} f_{n}(j)\leq \int_{A_n}^{\infty} \frac{(\log t)^n }{t^2}dt.$$
We now estimate 
$$I_n:=\int_{A_n}^{\infty} \frac{(\log t)^n }{t^2}dt=\int_{\log A_n}^\infty x^n e^{-x}dx\leq \int_{0}^\infty x^n e^{-x}dx=n!\leq n^n.$$
So that, simplifying by $n^n$: 
$$\mathbb{P}(\Omega_{n}(A_n))\lesssim 4(Bc)^n\leq 4\times 6^{-n},$$
and $\sum_{n\geq 1} \mathbb{P}(\Omega_{n}(A_n))\leq 4/5<1$ 
which proves (1). 

 For (2), by  Bernstein's deviation inequality and Borel-Cantelli, we know that almost surely (say for $\omega\in \Omega_1$), we have for  $n\geq n_{0}(\omega)$:
$$|R(\omega)\cap [1,A_n]|\leq 2\mathbb{E}(|R_{A_n}(\omega)|)=2m_{A_n}\leq 2c\log A_n \leq n. $$

This implies that, with $\Omega_2=\big(\Omega\setminus{\cup_{n\geq 1}} \Omega_{n}(A_n)\big)\cap \Omega_1$,  we have $\mathbb{P}(\Omega_2)>0$ with moreover: if $\omega\in \Omega_2$, then, for all $n$, we have 
\begin{equation}\label{fibre}R(\omega)\cap[A_n,\infty[ \hbox{\ contains no relation of length}\  \leq n.\end{equation}
 Otherwise, $R(\omega)\cap[A_n,\infty[$ would contain a relation of length $s$ with $3\leq s\leq n$ and since $A_n\geq A_s$, we would have 
$\omega\in \Omega_{s}(A_s)$, contradicting $\omega\in \Omega_2$.

 And moreover, for large $n$, say $n\geq n_{0}(\omega)$,  
\begin{equation}\label{tech} |R(\omega)\cap [1,A_n]|\leq n. \end{equation}
 Let now $E$ be a finite subset of $R(\omega)$ with  cardinality 
$|E|=2n$ or $2n+1$ such that 
$n\geq n_{0}(\omega)$. 
By (\ref{tech}) above, we know that 
$$|E\cap [1,A_n]|\leq |R(\omega)\cap [1,A_n]|\leq n. $$
So that 
$$|E\cap [A_n,\infty[| \geq |E|-|E\cap [1,A_n]|\geq 2n-n=n.$$
But if we now take 
$$F\subset E\cap [A_n,\infty[\subset R(\omega)\cap [A_n,\infty[$$
 with cardinality $n$, we see by definition (cf.~\ref{fibre}) that $F$ is quasi-independent. Since $|F|\geq (1/3)|E|$ and $E$ is arbitrary, this means, by Pisier's criterion, that $R(\omega)$ is Sidon. All in all, we proved that $R(\omega)$ is Sidon with positive probability. Since being Sidon is clearly an asymptotic property for $R(\omega)$, the zero-one law shows that $R(\omega)$ is Sidon almost surely.\\
 
 \noindent {\bf Step 2.} Dropping the dependence in $\omega$,  $R$ is almost surely a finite union of sets $R_j$ satisfying the assumptions of Step 1. Since a finite union of Sidon sets is again a Sidon set (Drury's theorem), $R$ itself is Sidon and we are done.
 
 For that, we select a large integer $M$ such that $C/M\leq 1/(48 e)$ and set 
 $$R_{j}=\{ \xi_{kM +j}: k\geq 1 \},\  j=0,1, \dots, M-1.$$
 For each fixed $j$, $R_j$ consists of  selectors of mean $\delta_{kM+j}$.  Clearly, since the $\delta_k$'s decrease,  for each $0\leq j\leq M-1$ we have
 $$\sum_{k=1}^N \delta_{kM} \leq \sum_{k=1}^N \delta_{kM-j} .$$
 Adding those inequalities gives 
 $$M\sum_{k=1}^N \delta_{kM}\leq \sum_{j=0}^{M-1}  \sum_{k=1}^N \delta_{kM-j}\leq \sum_{l=1}^{NM} \delta_l$$
 so that
 $$\sum_{k=1}^N  \mathbb{E}(\xi_{kM})= \sum_{k=1}^N \delta_{kM}\leq \frac{1}{M}  \sum_{l=1}^{NM} \delta_{l}\leq \frac{C}{M} \log (NM)\leq  \frac{2C}{M} \log N$$
 for  $N\ge M$.  As $2C/M\leq 1/(24e)$,
    $R_M$ satisfies the assumptions of Step 1, and is almost surely Sidon.  We do the same for $R_j,  \ j=0,1,\dots, M-1$ and we are done. 
    \end{proof} 
    \smallskip
      \noindent   {\bf Remark.}
      \smallskip
    \noindent   Even for $(\delta_n)$ nonincreasing, the assumption $m_N=O(\log N)$  
    is more general than the Kahane-Katznelson assumption $\delta_n=O(1/n)$, as shown by the following example.  First we 
    choose $x_j=2^{2^j}$ so that 
   $$x_{j+1}=x_{j}^{2}, \ \ \log x_j\asymp 2^ j,\ \  x_1=4.$$ 
 Then we define 
  \begin{equation*}
    \delta_n=
    \begin{cases} 
    1 \text{\quad if} \quad 1\le n\leq 4,\\ 
    \frac{\log x_{j+1}}{{x_{j +1}}} \text{\quad if}\quad  x_j<n\leq x_{j+1} \hbox{\ with}\ j\geq 2.
    \end{cases} 
     \end{equation*}
     Observe that\\
   (i)  The sequence $(\delta_n)$ is clearly \textit{nonincreasing.}\\
   (ii) The assumption  $\delta_n=O(1/n)$  fails  for $n=x_{j+1}$. \\
   (iii) Finally, for large $N$, let $n$ satisfy $x_n<N\leq x_{n+1}$. Then
   $$ m_N= \sum_{k=1}^N \delta_k \lesssim \sum_{j=1}^n \Big(x_{j+1} \frac{\log x_{j+1}}{x_{j+1}}\Big)= \sum_{j=1}^n  \log x_{j+1}
    \lesssim \sum_{j=1}^n 2^{ j} \lesssim  \log N.$$

\section{Questions on $S$ and $T$}

 \begin{Pb}  Theorem \ref{S_N} and Theorem \ref{ss} can be stated and proved in the same way for $S(q_1, q_2)$. But efforts  are needed  for $S(q_1, q_2, q_3)$.
\end{Pb}

\begin{Pb} Is $S$ a rigid set ?
\end{Pb}
\begin{Pb} Is a Sidon set always rigid ?
\end{Pb}

\begin{Pb}    Is a Bohr-dense set always uniformly recurrent ?                   
\end{Pb}

\begin{Pb} Is there a Bohr-dense Sidon set ? Equivalently, a non Bohr-closed Sidon set?
\end{Pb}

\begin{Pb} Is the Furstenberg set $S$ $\frac{4}{3}$-Rider (even $\frac{4}{3}$-Sidon) ?
More generally,  is $S(q_1,\ldots, q_s)$  $\frac{2s}{s+1}$-Rider (even   $\frac{2s}{s+1}$-Sidon) ?
 
\end{Pb}

  \begin{Pb}
    Is $T$  a $\frac{4}{3}$-Rider set ?
  \end{Pb}

\begin{Pb}
 A simple argument  
   shows that, in Theorem \ref{mieux que kk}, we can replace $m_N=O(\log N)$ by $m_{N_j}=O(\log N_j)$ where $(N_j)$ is an increasing sequence of integers such that $N_{j+1}=O(N_j)$. 
   But we do not know if 
   the assumption 
     $\liminf_{N\to \infty}{m_N}/{\log N}<\infty$ is enough.

\end{Pb}

\subsection*{Acknowledgments.} 
A. H. Fan is partially supported by NSFC (grant no.11971192 and grant no. 12231013). 
 H.~Queff\'elec and M.~Queff\'elec acknowledge the support of the Labex CEMPI (ANR-11-LABX-0007-01).



\begin{thebibliography}{99}


\bibitem{Ar} P. Arnoux, {\it Sturmian sequences}. Substitutions in dynamics, arithmetics and combinatorics, 143--198, Lecture Notes in Math., 1794, Springer, Berlin, 2002. 

\bibitem{BG} C. Badea, S. Grivaux, \textit{Kazhdan constants, continuous probability measures with large Fourier coefficients and rigidity sequences,} Comment. Math. Helv. 95 (2020), no. 1, 99-127.

\bibitem{BF2005} J. Barral and A. H. Fan 
  \textit{Covering numbers of different points in {D}voretzky covering}, 
    Bull. Sci. Math.,
  vol. 129, no. 4, (2005), 275-317. 


\bibitem{BDH}  A. B\'erczes, A. Dujella, L. Hajdu, \textit{Some Diophantine properties of the sequence of $S$--units.} J. Number Theory 138 (2014), 48--68.
\bibitem{BJLR} V. Bergelson, A. del Junco, M. Lemanczyk, J. Rosenblatt,  \textit{Rigidity and non-recurrence along sequences. }
Ergodic Theory Dynam. Systems 34 (2014), 1464--1502. 
\bibitem{BL} V. Bergelson, E. Lesigne, \textit{ Van der Corput sets in $\Z^d$}. Colloq. Math. 110 (2008),  1-49.

\bibitem{Bewley1971} 
 T. Bewley, \textit{ Extension of the Birkhoff and von Neumann ergodic theorems to semigroup actions}. Ann. Inst. H. Poincaré Sect. B (N.S.) 7 (1971), 283-291.

\bibitem{Bl} R. Blei, \textit{Analysis in Integers and Fractional dimensions,} Cambridge 71, 2001. 
Ark. Mat., 36 (1998), 307-316.


\bibitem{BLM} S. Boucheron, G.~Lugosi, P.~Massart, \textit{On concentration of self-bounding functions,} Electronic J.~Probab. 14 (64)  (2009),  1884-1899.

\bibitem{Bo1} J. Bourgain, Ruzsa's problem on sets of recurrence, Isr. J. Math. 59 (1987), 150--166.

\bibitem{Bo2} J. Bourgain,  \textit{On the maximal ergodic theorem for certain subsets of the integers,} Isra\"el J.~Math.~61 (1988), 39-72.
\bibitem{Bug 2} Y. Bugeaud,  \textit{Distribution modulo one and diophantine approximation}, Cambridge University Press, Vol.193, 2012. 
\bibitem{Bu} D.~Burkholder,\,\textit{Martingale transforms}, Ann.~Math.~Stat., {37} (1966), 1494-1504.
\bibitem{Bu 2} D.~Burkholder,\,\textit{Sharp inequalities for martingales}, Ast\'erisque 157-158 (1988),  75-94.

\bibitem{CF}  Ch. Cuny and A. H. Fan, \textit{Study of almost everywhere convergence of series by mean of martingale methods}. 
Stochastic Process. Appl. 127 (2017), no. 8, 2725-2750. 

\bibitem{Dr} S. Drury,  \textit{Sur les ensembles de Sidon,} C.~R.~Acad.~Sc.Paris {t. 271} (1970), 162-3.
\bibitem{EG} T. Eisner, S. Grivaux, \textit{Hilbertian Jamison sequences and rigid dynamical systems.} 
J. Funct. Anal. 261 (2011), 2013--2052. 



\bibitem{ET}  P. Erd\"os,  S. J. Taylor, \textit{On the set of points of convergence of a lacunary trigonometric series and the equidistribution properties of related sequences}, Proc. London Math. Soc. 7,  598-615, 1957. 

\bibitem{Fal} K. Falconer,  \textit{Fractal Geometry,}
John Wiley and Sons (1990).

\bibitem{Fa} A. H. Fan,  \textit{Lacunarit\'e \`a la Hadamard et \'equir\'epartition,}
Colloq. Math. 66 (1993), 151-163.



\bibitem{AHM}  A.~H.~ Fan, H.~Queff\'elec, M.~Queff\'elec,
\textit{The Furstenberg set and its random version,}
L'Enseignement Math\'ematique, to appear.\\
https://ems.press/journals/lem/articles/8188703 (online first)

\bibitem{Fl} L.~Flaminio,  \textit{Private communication,}
University of Lille, December 2022.

\bibitem{FS} A. H. Fan, D. Schneider,  \textit{Recurrence properties of sequences  of integers,} Sci. China Math. 53 (2010) no.3, 641-656. 








\bibitem{FT} E. Fouvry, G. Tenenbaum,  \textit{Entiers sans grand facteur premier en progressions arithm\'etiques,}  Proc. London Math. Soc. 63 (1991), no. 3, 449--494.

\bibitem{Fu1} H. Furstenberg, \textit{Disjointness in ergodic theory, minimal sets, and a problem in Diophantine approximation,} Math. Systems Theory 1 (1967), 1-49.

\bibitem{Ge} A.~Gelfond, \textit{M\'ethodes \'el\'ementaires dans la th\'eorie analytique des nombres,} Gauthier-Villars 1965. 

\bibitem{Gr} J.~Griesmer, \textit{Separatig Bohr's denseness from measurable recurrence,} Discrete Anal. 2021, Paper No 9, 2O pp.

\bibitem{GV} R.~Gundy, N.~Varopoulos, \textit{A martingale that occurs in harmonic analysis}, Arkiv.~Math.~14 (1976) ({2}), 179-187.

\bibitem{HaWr} G.~Hardy, E.~Wright, \textit{An introduction to the theory of numbers},  Oxford, 1968.

\bibitem{Ha} G.~H.~Hardy, \emph{Ramanujan, Twelve lectures on subjects suggested by his life and work}, Cambridge, 1940.
\bibitem{Har} S. Hartman, {\it The method of Grothendieck-Ramirez and weak topologies in C(T),} 
Studia Math. 44 (1972), 181--197.

\bibitem{Hu} Y. Huang, {\it Ergodic theorems for random sets with density zero,} 
Ergod. Th. - Dynam.~Sys. (1994), {\bf 14}, 141--149.

\bibitem{Hof} J. Hoffmann-Jorgensen, {\it Sums of independent Banach-space valued random variables,} 
Studia Math. 44 (1974), 159--186.



\bibitem{JoSc}B.~Johnson, G.~Schechtman,  \textit{Remarks on Talagrand's deviation inequality for Rademacher functions,} Lecture Notes in Math., Longhorn Notes, Springer, Berlin, 1991.

\bibitem{JW} G.~Johnson, G.~Woodward, \textit{On $p$-Sidon sets,} Indiana Univ. Math. J. { 24},  161-167, 1974.
\bibitem{Ka} J. P.~Kahane, \textit{Some random series of functions},  Cambridge 1985.
\bibitem{KK1} J. P.~Kahane, Y.~Katznelson, \textit{Entiers al\'eatoires et analyse harmonique,} J. Anal.~Math.~ 105 (2008), 363-378.
\bibitem{KK2} J. P.~Kahane, Y.~Katznelson, \textit{Distribution uniforme de certaines suites d'entiers al\'eatoires dans le groupe de Bohr,} J. Anal.~Math. 105 (2008), 379-382.
\bibitem{KS} J. P.~Kahane, R.~Salem, \textit{Ensembles parfaits et s\'eries trigonom\'etriques,} Second Edition,  Hermann,  1994.
\bibitem{Kat} Y. Katznelson, {\it Sequences of integers dense in the Bohr group.} Proc. Roy. Inst. of Tech. (1973), 79--86.
\bibitem{Kh} A. Khinchin, \textit{Ein Satz \"uber Kettenbruche mit arithmetischen Anwendungen,} Math.~Zeit. 18 (1923), 289-306.
\bibitem{Ko} J.~F.~Koksma,  \textit{A diophantine property  of summable functions,}  Indian journal of Mathematics Society {15} (1951), 87-96.
\bibitem{KN} L.~Kuipers, H.~Niederreiter, \textit{Uniform distribution of sequences,} Dover publications, Second edition, 2002.


\bibitem{LR} P.~Lef\`evre, L.~Rodriguez-Piazza, \textit{$p$-Rider sets are $q$-Sidon sets,} PAMS {131} (6) (2003), 1829-1838.

\bibitem{LoR} J.~Lopez, K.~Ross, \textit{Sidon sets.} Lecture Notes in Pure and Applied Mathematics, Vol.13, 1975.

\bibitem{LQ} D.~Li, H.~Queff\'elec, \textit{Introduction to Banach spaces: Analysis and Probability,} Cambridge Studies in Advanced Mathematics (2018).
\bibitem{LQR} D.~Li, H.~Queff\'elec, L.~Rodr\'iguez-Piazza, \textit{Some new thin sets of integers in harmonic analysis,} Journal d'Analyse Math\'ematique {86} (2002), 105-138. 

\bibitem{LQR2} D.~Li, H.~Queff\'elec, L.~Rodr\'iguez-Piazza, \textit{On some random thin sets of integers,} PAMS  {136} (2008), 141-150. 



\bibitem{Ly0}  R.~Lyons,  oral communication.
\bibitem{Ly1} R. Lyons, \textit{On measures simultaneously 2- and 3-invariant.} Israel J. Math. 61 (1988), 219--224. 


\bibitem{Mane} R.~Ma\~{n}\'e,  \textit{ Ergodic theory and differentiable dynamics,} Springer-Verlag (1987). 

\bibitem{Ma} J.~M.~Marstrand,  \textit{ On Khinchin's conjecture about strong uniform distribution,} Proc.~London Math.~Soc. 21 (1970), 540-556. 

\bibitem{MS} C.~Muscalu and W.~Schlag, \emph{Classical and multilinear harmonic Analysis, Volume I,} Cambridge studies in advanced mathematics 137,  2013.


\bibitem{Na} R.~Nair, \textit{On strong uniform distribution,} Acta Arithmetica LVI (1990), 183-192.



\bibitem{Ne2} J. Neveu, \textit{Martingales \`a temps discret,} Masson 1971.

\bibitem{Pa} R.C.~ Paley, \textit{On the lacunary coefficients of power series,} Ann. of Math. (2) 34 (1933), no 3, 615-616.


 
\bibitem{Ph} W.~Philipp, \textit{Empirical distribution functions and strong approximation theorems for dependent random variables. A problem of Baker in probabilistic number theory}, Trans.
Amer. Math. Soc. 345 (1994), 705-727. 


\bibitem{Pi} G.~Pisier, \textit{De nouvelles caract\'erisations des ensembles de Sidon,} Math. Anal. and Applic., Part B, Advances in Math.~Suppl.~Studies, Vol {7 B}, 1981. 

\bibitem{Pi2} G.~Pisier, \textit{Martingales in Banach spaces,} Cambridge University Press, 2016.



 \bibitem{Re}A.~R\'enyi,  \textit{Calcul des Probabilit\'es,} Dunod, 1966. 


 \bibitem{Rh}G.~Rhin,  \textit{Approximants de Pad\'e et mesures effectives d'irrationalit\'e,} S\'eminaire de th\'eorie des nombres, Paris 1985-1986, 155-164. Progr.~Math.~71, Birkha\"user Boston, Boston, MA, 1987.
 
\bibitem{Ri}D.~Rider,  \textit{Randomly continuous functions and Sidon sets,} Duke Math.~J.~ { 42} (1975), 759-764.

\bibitem{Ro1} L. Rodriguez-Piazza, \textit{Caract\'erisation des ensembles $p$-Sidon p.s.}, C.~R.A.~S. Paris S\'er.~I Math.~305 (1987), 237-240.

\bibitem{Ro2} L. Rodriguez-Piazza, \textit{Rango y propriedades de medidas vectoriales: Conjuntos $p$-Sidon p.s.}, Thesis, Universidad de Sevilla, 1991.

\bibitem{Ros} J.~Rosenblatt, \textit{Universally bad sequences in ergodic theory}, in: Almost everywhere convergence II. Proc.2nd Int.Conf.,Evanston /IL (USA) $1989$, Academic Press, Boston, MA, 1991, pp.227--245.

\bibitem{Ru1}W.~Rudin,  \textit{Fourier Analysis on groups,}  second edition, Cambridge 1993.

\bibitem{Ru2}W.~Rudin,  \textit{Trigonometric series with gaps,} J.~Math.~Mech.~{9}, 203-227, 1960.
\bibitem{Rud}D. J. Rudolph, {\it $\times 2$ and $\times 3$ invariant measures and entropy,} Ergodic Theory Dynam. Systems 10 (1990), no. 2, 395--406. 
\bibitem{Ru} I.~Z. Ruzsa, Connections between the uniform distribution of a sequence and its difference, Coll. Math. Soc. Janos Bolyai, Topics in classical number theory, Budapest, 1981.



\bibitem{Tal}M.~Talagrand,  \textit{An isoperimetric theorem for the cube and the Khintchine-Kahane inequalities,}Proc.~Amer.~Math.~Soc.~ Vol.~ {104}(3) (1988), 1905--910.

\bibitem{Ti1} R. Tijdeman, \textit{On integers with many small prime factors,} Compos. Math. 26 (1973), 319--330.
\bibitem{Ti2} R. Tijdeman, \textit{On the maximal distance between integers composed of small primes,} Compos. Math. 28 (1974), 159--162.

\bibitem{Us} S.~Usuki, \textit{$\times a\times b$ empirical measures, the irregular set and entropy,} arXiv:2205. 06605 v 2 21 Jul 2022. 

\bibitem{WaWu} L.~Wang, Q.~Wu, \textit{On the irrationality measure of $\log 3$, } J.~Number Theory 142 (2014), 264-273. 

\bibitem{Wo} G. Woodward, \textit{$p$-Sidon sets and a uniform property}, Indiana Univ.~Math.~J. {25} (6), p.~995-2003, 1976.
\bibitem{Zy} A.~Zygmund, \textit{Trigonometric series,} sixth edition, Cambridge 1993.

\end{thebibliography}
\end{document}